\documentclass[ejs,preprint,dvips,twoside,linksfromyear]{imsart}

\RequirePackage[OT1]{fontenc}
\RequirePackage{amsthm,amsmath,amssymb}
\RequirePackage[numbers]{natbib}
\RequirePackage[colorlinks,citecolor=blue,urlcolor=blue]{hyperref}

\startlocaldefs
\numberwithin{equation}{section}
\theoremstyle{plain}
\newtheorem{thm}{Theorem}[section]

\newtheorem{lemma}{Lemma}[section]
\newtheorem{remark}{Remark}[section]

\newtheorem{proposition}{Proposition}[section]
\endlocaldefs

\begin{document}
\begin{frontmatter}
\title{On maxisets of nonparametric tests\thanksref{T1}}
\thankstext{T1}{Supported in part RFFI Grant 17-01-00828}
\runtitle{On maxisets of nonparametric tests}

\begin{aug}
\author{\fnms{Mikhail} \snm{Ermakov}\ead[label=e1]{erm2512@gmail.com}}

\address{Institute of Problems of Mechanical Engineering RAS, Bolshoy pr., 61, VO, 1991178 RUSSIA and
St. Petersburg State University, Universitetsky pr., 28, Petrodvoretz, 198504 St. Petersburg, RUSSIA\\
\printead{e1}}
\runauthor{M. Ermakov}
\end{aug}
\begin{abstract} For the problems of nonparametric hypothesis testing we introduce  the notion of maxisets and maxispace. We point out the maxisets  of $\chi^2-$tests, Cramer-von Mises tests, tests generated  $\mathbb{L}_2$- norms of kernel estimators and  tests generated quadratic forms of estimators of Fourier coefficients. For these tests we show  that, if sequence of alternatives having given rates of convergence to hypothesis is consistent, then  each altehrnative can be broken down into the sum of two parts: a function belonging  to maxiset and orthogonal function. Sequence of functions belonging to maxiset is consistent sequence of alternatives.
  We point out asymptotically minimax tests if sets of alternatives are maxiset with deleted "small" $\mathbb{L}_2$-balls.
\end{abstract}

\begin{keyword}[class=AMS]
\kwd[Primary]\,{62F03} \kwd{62G10}  \kwd{62G20}

\end{keyword}

\begin{keyword}
\kwd{Maxisets}
\kwd{ chi-squared test}
\kwd{consistency}
\kwd{nonparametric hypothesis testing}
\kwd{signal detection}
\end{keyword}

\end{frontmatter}
\section{Introduction}

Let $X_1,\ldots,X_n$ be i.i.d.r.v.'s with c.d.f. $F(x)$, $x \in (0,1)$. Let c.d.f. $F(x)$ have a density $p(x) = 1 + f(x) = dF(x)/dx, x \in (0,1)$. Suppose that $f \in \mathbb{L}_2(0,1)$ with the norm
$$
\|f\| = \left(\int_0^1 f^2(x) dx \right)^{1/2} < \infty.
$$
We explore the problem of testing hypothesis
\begin{equation}\label{i1}
{\rm H_0}: f(x) = 0, \quad x \in (0,1)
\end{equation}
versus nonparametric alternatives
\begin{equation}\label{i2b}
{\rm H_n} : f \in V_n = \{ g : \|g\| \ge cn^{-r},\, g \in U\},
\end{equation}
where $U$ is a ball in some functional space $\Im \subset \mathbb{L}_2(0,1)$. Here $c,r$ are constants, $c>0,  0< r <1/2$.

We could not verify the hypothesis  ${\rm H_0}$ if  nonparametric sets of alternatives contains all densities $p$, $\|p-1\| \ge cn^{-r}$, (see Le Cam and Schwartz \cite{les}, Ingster \cite{ing02}, Ermakov \cite{er15}). The problem can be solved if
additional a priori information  is provided that the function $f$ belongs to a compact ball $U$.  For the problems of hypothesis testing in functional spaces the surveys of the results exploring  the setup  (\ref{i1}) and (\ref{i2b}) one can find in   Horowitz and Spokoiny \cite{hor}, Ingster and Suslina \cite{ing02},
 Laurent, Loubes and Marteau \cite{la} and
 Comminges and Dalalyan \cite{dal} (see also references therein). 
 
 The problem of asymptotically minimax nonparametric estimation is also  explored if a priori information is provided that unknown function belongs to some set $U$. In this setup the set $U$ is a  compact in some functional space (see Le Cam and Schwartz \cite{les},  Ibragimov and Khasminskii \cite{ih}, Johnstone \cite{jo}).

 The paper goal is to find the largest functional spaces $\Im$ allowing to test these hypotheses if we  implement one  of widespread test statistics presented below.   The  largest space $\Im$ we call maxispace. The balls $U \subset \Im$ with center at zero we   call maxisets.

There are few results related to the study of rate of consistency of traditional nonparametric tests if the sets of alternatives are nonparametric. First of all we should mention Mann and Wald paper \cite{mann}. Mann and Wald \cite{mann} established the optimal order of number of cells for chi-squared tests if Kolmogorov distances of alternatives to hypothesis are greater some constants. If   $U$ is a ball in Besov space $\mathbb{B}^s_{2\infty}$, the problem of testing hypothesis $H_0$ versus alternative $H_n$ has been explored   Ingster \cite{ing87} for chi-squared tests with increasing number of cells, Kolmogorov and Cramer - von Mises tests. Horowitz and Spokoiny \cite{hor}  and Ermakov \cite{er97, er03, er11} explored asymptotically minimax properties of wide-spread nonparametric tests in semiparametric setup based on the distance method.

 In paper we show that  Besov spaces $\mathbb{B}^s_{2\infty}$ are maxispaces for $\chi^2$-tests,  Cramer- von Mises tests  and tests generated  $\mathbb{L}_2$- norms of kernel estimators. For the problem of signal detection in Gaussian white noise, for tests generated quadratic forms of estimators of Fourier coefficients, we show that the assignment of maxispaces in  orthonormal basis coincides with the assignment  of  Besov spaces $\mathbb{B}^s_{2\infty}$ in trigonometric basis.

 The tests generated  $\mathbb{L}_2$-- norms of kernel estimators and  the tests generated quadratic forms of estimators of Fourier coefficients are explored for the problem of signal detection in Gaussian white noise.
  We observe a realization of random process $Y_n(t)$ defined stochastic differential equation
\begin{equation}\label{q1}
dY_n(t) = f(t) dt + \frac{\sigma}{\sqrt{n}}\, dw(t), \quad t \in [0,1],\quad \sigma >0,
\end{equation}
where  $f \in \mathbb{L}_2(0,1)$  is  unknown signal and $dw(t)$ is Gaussian white noise.

 The problem of hypothesis testing is the same.

  This setup allows  do not make additional assumptions and  to simplify the reasoning.  More traditional  problems of hypothesis testing are explored for i.i.d.r.v.'s.

  For  nonparametric estimation the notion of maxisets has been introduced Kerkyacharian and Picard \cite{ker93}. The maxisets of widespread nonparametric  estimators have been comprehensively explored (see Cohen,  DeVore, Kerkyacharian,   Picard \cite{co},  Kerkyacharian and Picard \cite{ker02}, Rivoirard \cite{rio},  Bertin and Rivoirard \cite{rio09}, Ermakov \cite{er17} and references therein).

  $\mathbb{L}_2$--norm  is naturally arises in the study of test behaviour  for the problems of hypothesis testing with  alternatives converging to hypothesis.
If we consider the problem of testing hypothesis (\ref{i1}) versus simple alternatives $H_{1n}: f(x) =f_n(x) =  c\, n^{-1/2} h(x)$, $\|h\|<\infty$, then the asymptotic of type II error probabilities of Neymann-Pearson tests is defined by $\|h\|^2$. Similar situation takes place  also for the problem of signal detection in Gaussian white noise.

For the study of maxiset properties we introduce the notion of perfect maxisets. The definition of perfect maxiset is based on the notion of $n^{-r}$--consistency (see Tsybakov \cite{ts}).
Consistent sequence of alternatives having $n^{-r}$--rate of convergence  to hypothesis is called $n^{-r}$-consistent.

The maxisets are called perfect if the functions of any $n^{-r}$-consistent sequence of alternatives admits the representation as the sum of  two functions:   functions from some maxiset and the orthogonal functions. The sequence of alternatives corresponding to the functions from maxiset is $n^{-r}$-consistent. Moreover the sum of functions of $n^{-r}$-consistent sequence of alternatives belonging to maxiset and functions of sequence of inconsistent alternatives having $n^{-r}$-rate of convergence to hypothesis  also form $n^{-r}$-consistent sequence.

Thus  all information on $n^{-r}$--consistency of  sequences of alternatives   is contained in perfect maxisets.
We show that the  maxisets of all above mentioned tests are perfect.

We introduce also the notion of completely perfect maxisets.  If the notion of perfect maxisets requires the exploration of type I and type II error probabilities in terms of rates of convergence, the exploration of completely perfect maxisets requires the knowledge of strong asymptotic of type I and type II error probabilities. We show that maxisets of Cramer--  von Mises tests, tests generated  $L_2$-- norms of kernel estimators and  tests generated quadratic forms of estimators of Fourier coefficients are completely perfect.

Another nonasymptotic approach to the maxiset definition has been proposed recently Autin, F., Clausel,M., Jean-Marc Freyermuth, J. and  Marteau \cite{au}.

Paper is organized as follows. In section 2 we discuss desirable properties of maxisets and maxispaces. After that we provide the definitions of maxisets and maxispace, perfect and completely perfect maxisets. In section 3, under rather general assumptions,  we show that, if there is comsistent test, then the set $U$ should be compact. In sections  4, 5, 6 and 7 we point out maxisets of  test statistics based on quadratic forms of estimators of Fourier coefficients,  $\mathbb{L}_2$ -- norms of kernel estimators, $\chi^2$--tests and Cramer-- von Mises tests   respectively. In these sections we explore also the properties of these maxisets.  In section 8 we point out asymptotically minimax test statistics, if a priori information is provided, that  alternatives  belong to maxiset with "small $\mathbb{L}_2$ -balls removed". Sections 9  contains the proofs of all Theorems.

We use letters $c$ and $C$ as a generic notation for positive constants. Denote $\chi(A)$ the
indicator of an event $A$.  Denote $[a]$ the whole part of real number $a$. For any two sequences of positive real numbers $a_n$ and $b_n$, $a_n = O(b_n)$ and $a_n \asymp b_n$ imply respectively $a_n < C b_n$ and $ca_n \le b_n \le Ca_n$ for all $n$.

Denote
$$ \Phi(x) = \frac{1}{\sqrt{2\pi}}\,\int_{-\infty}^x\,\exp\{-t^2/2\}\, dt, \quad x \in \mathbb{R}^1,
$$
the standard normal distribution function.

Let $\phi_j, 1 \le j < \infty$, be orthonormal system of functions. Define the sets
\begin{equation}\label{vv}
\mathbb{\bar B}^s_{2\infty}(P_0) = \Bigl\{f : f = \sum_{j=1}^\infty\theta_j\phi_j,\,\,\,  \sup_{\lambda>0} \lambda^{2s} \sum_{j>\lambda} \theta_j^2 < P_0,\,\, \theta_j \in \mathbb{R}^1 \Bigr\}. 
\end{equation}
Under some conditions on the basis $\phi_j, 1 \le j < \infty,$  the space
$$
\bar{\mathbb{ B}}^s_{2\infty} = \Bigl\{ f : f = \sum_{j=1}^\infty\theta_j\phi_j,\,\,\,  \sup_{\lambda>0} \lambda^{2s} \sum_{j>\lambda}\, \theta_j^2 < \infty,\,\, \theta_j \in \mathbb{R}^1 \Bigr\}.
$$
is Besov space $\mathbb{B}^s_{2\infty}$ (see  Rivoirard  \cite{rio}).
In particular, $\mathbb{\bar B}^s_{2\infty}$ is Besov space  $\mathbb{B}^s_{2\infty}$if $\phi_j$, $1 \le j < \infty$, is trigonometric basis.

If $\phi_j(t) = \exp\{2\pi i j t\}$, $x\in (0,1)$, $j = 0, \pm 1, \ldots$, is trigonometric basis, denote
$$
\mathbb{ B}^s_{2\infty}(P_0) = \Bigl\{f : f = \sum_{j=-\infty}^\infty \theta_j\phi_j,\,\,\,  \sup_{\lambda>0} \lambda^{2s} \sum_{|j| >\lambda} |\theta_j|^2 < P_0 \Bigr\}.
$$
The balls in Nikols'ki classes
$$
\int\,(f^{(l)}(x+t) - f^{(l)}(x))^2\, dx \le L |t|^{2(s-l)}, \quad \|f\| < C
$$
with $l = [s]$ are the  balls in $\mathbb{B}^s_{2\infty}$.

We also introduce definition of balls in Besov spaces  $\mathbb{B}^s_{2\infty}$ in terms of wavelet basis $\phi_{kj}(x) = 2^{(k-1)/2}\phi(2^{k-1}x-j),\, 1 \le j < 2^k,\, 1 \le k < \infty$.
Denote
$$
\mathbb{\tilde B}^s_{2\infty}(P_0) = \Bigl\{f : f = 1 + \sum_{k=1}^\infty\sum_{j=1}^{2^k}\,\theta_{kj}\,\phi_{kj},\,\,\,  \sup_{\lambda>0}\,2^{2\lambda s} \sum_{k>\lambda}\,\sum_{j=1}^{2^k}\,\theta^2_{kj} \le P_0,\,\theta_{jk} \in \mathbb{R}^1 \Bigr\}.
$$
\section{Main definitions}
\subsection{Definition of consistency and $n^{-r}$-consistency}
For any test $K_n= K_n(X_1,\ldots,X_n)$  denote $\alpha(K_n)$ its type I error probability, and $\beta(K_n,f)$ its type II error probability for the alternative $f \in \mathbb{L}_2(0,1)$.

For the problem of testing hypothesis ${\rm H_0}\,:\, f=0$ versus alternatives ${\rm H_n}\, :\, f=  f_n$, we say that sequence of alternatives $f_n$ is consistent  if there is sequence of tests $K_n$ generated test statistics $T_n$ such that
\begin{equation}\label{}
\limsup_{n\to\infty} (\alpha(K_n) + \beta(K_n, f_n)) < 1.
\end{equation}
If $cn^{-r} < \|f_n\| < Cn^{-r}$ additionally, we say that sequence of alternatives $f_n$ is $n^{-r}$- consistent.

We say that sequence of alternatives $f_n$ is inconsistent  if for each sequence of tests $K_n$ generated test statistics $T_n$ there holds
\begin{equation}\label{}
\limsup_{n\to\infty} (\alpha(K_n) + \beta(K_n, f_n)) \ge 1.
\end{equation}
If $cn^{-r} < \|f_n\| < Cn^{-r}$ additionally, we say that sequence of alternatives $f_n$ is $n^{-r}$- inconsistent.

Denote
$$
\beta(K_n,V_n) = \sup\{ \beta(K_n,f), f \in V_n\}.
$$
We say that, for test statistics $T_n$, the problem of hypothesis testing is $n^{-r}$-consistent on  the set $U$ (consistent on the sets $V_n$ respectively) if there is sequence of tests $K_n$ generated test statistics $T_n$ such that
\begin{equation}\label{}
\limsup_{n\to\infty} (\alpha(K_n) + \beta(K_n,V_n)) < 1.
\end{equation}
\subsection{Definition of maxisets and maxispace}
Let us discuss desirable properties of maxisets and maxispaces of test statistics $T_n$ based on similar definition of maxisets in estimation (see Kerkyacharian and Picard \cite{ker93, ker02}).

We would like  to find  Banach space $\Im\subset \mathbb{L}_2(0,1)$ such that
\vskip 0.3cm
{\sl i.} problem of hypothesis testing is $n^{-r}$-consistent on the balls  $U \subset \Im$
\vskip 0.3cm
{\sl ii}. for any $f \notin \Im, f \in \mathbb{L}_2(0,1)$, for tests $K_n, \alpha(K_n)  = \alpha(1+o(1))$, $0  < \alpha < 1$, generated test statistics $T_n$, there are  functions $f_{1n},\ldots,f_{k_nn} \in \Im$  such that
$$
 cn^{-r}\le \Bigl\|f - \sum_{i=1}^{k_n} f_{in}\Bigr\| \le Cn^{-r}
$$
and
\begin{equation}\label{er1}
\limsup_{n \to \infty}\beta\Bigl(K_n,f - \sum_{i=1}^{k_n} f_{in}\Bigr) \ge 1-\alpha,
\end{equation}
{\sl iii}. the space $\Im$ contains  smooth functions up to the functions of "the smallest possible smoothness" for this setup.
\vskip 0.3cm
Let us discuss the content of the second point {\sl ii.} of this definition.
We could not proof such a statement for arbitrary functions $f_{in} \in \Im$.  We shall suppose that functions $f_{in}$  belong to specially defined finite dimensional subspaces $\Pi_k \subset \Im$.

Let us discuss the third point  {\sl iii.} of desirable definition.
 We can take arbitrary finite number  of unsmooth functions and search for the maxispace $\Im$ containing these functions.
Thus the maxispace problem is ambiguously defined without {\sl iii.}

The definition of maxisets  we begin with  preliminary notation.

Let $\Im \subset \mathbb{L}_2(0,1)$ be Banach space with  norm $\|\cdot\|_\Im$ and let $U(\gamma)=\{f: \|f\|_\Im \le \gamma, f \in \Im\}, \gamma > 0,$ be a ball in $\Im$.

Define subspaces $\Pi_k$, $1 \le k < \infty$, by induction.

Denote $d_1= \max\{\|f\|, f \in U(1)\}$ and denote $e_1$ function $e_1 \in U=U(1)$ such that $\|e_1\|= d_1.$ Denote $\Pi_1$ linear space generated vector $e_1$.

For $i=2,3,\ldots$ denote
$d_i = \max\{\rho(f,\Pi_{i-1}), f \in U \}$ with $\rho(f,\Pi_{i-1})=\min\{\|f-g\|, g \in \Pi_{i-1} \}$. Define function $e_i$, $e_i \in U$, such that $\rho(e_i,\Pi_{i-1}) = d_i$.
Denote $\Pi_i$ linear space generated functions $e_1,\ldots,e_i$.

For any $f \in \mathbb{L}_2(0,1)$ denote
$f_{\Pi_i}$ the projection of  $f$ onto the subspace $\Pi_i$ and denote $\tilde f_i = f - f_{\Pi_i}$.

Thus we associate with each $f \in \mathbb{L}_2(0,1)$ sequence of functions $\tilde f_i, \tilde f_i \to 0$ as $i \to \infty$. This allows to cover by  our consideration the all space $\mathbb{L}_2(0,1)$. Suppose that the functions $e_1,e_2,\ldots$ are sufficiently smooth. Then, considering the functions $\tilde  f_i = f - f_{\Pi_i}$, we "in some sense delete the most smooth part $f_{\Pi_i}$ of function $f$ and explore the behaviour of remaining part."

For the problem of hypothesis testing on a density we suppose  that for all $f   \in \Im$ there holds $\int_0^1 f(s) \, ds =0$.

We say that  $U(\gamma)$, $\gamma > 0$, is maxiset for test statistics $T_n$ and $\Im$ is maxispace if the following two statements take place
\vskip 0.3cm
{\sl i.} there is sequence of tests $K_n$, $\alpha(K_n) = \alpha(1+ o(1))$, $0<\alpha<1,$ generated test statistics $T_n$ satisfying the following inequality
\begin{equation}
\label{}
\limsup_{n\to\infty} (\alpha(K_n) + \beta(K_n,V_n)) < 1,
\end{equation}
\vskip 0.3cm
{\sl ii.}  for any  $f\in \mathbb{L}_2(0,1)$, $f \notin \Im$, $\int_0^1 f(s) \, ds =0$, there are sequences  $i_n, j_{i_n}$ with $i_n \to \infty$ as $n \to \infty$   such that $c j_{i_n}^{-r}<\|\tilde f_{i_n}\| < C j_{i_n}^{-r}$ for some constants $c$ and  $C$, and, if $1 + \tilde f_{i_n}(s) \ge 0$ for all $s \in [0,1]$, then any  sequence of tests $K_{j_{i_n}}$, $\alpha(K_{j_{i_n}}) = \alpha(1+ o(1))$, $0<\alpha<1,$ generated test statistics $T_{j_{i_n}}$  satisfies the following inequality
\begin{equation}
\label{}
\liminf_{n\to\infty} (\alpha(K_{j_{i_n}}) + \beta(K_{j_{i_n}},\tilde f_{i_n})) \ge 1.
\end{equation}
\vskip 0.3cm
All  definitions given above were provided in general terms. In each setup  these definitions are slightly different.
 For the problems of signal detection the requirements in definition of maxiset $1 + \tilde f_{i_n}(s) \ge 0$ for all $s \in [0,1]$ and $\int_0^1 f(s) \, ds =0$ are omitted. In definition of maxiset we replace indices $i=1,2,\ldots$ with $i = \pm 1, \pm 2, \ldots$ for the kernel- based tests and for the chi-squared tests.
\subsection{Definition of perfect maxisets and completely perfect maxisets}
We could not verify (\ref{er1}) for fixed $f \notin \Im$ and arbitrary $f_{in}  \in \Im$. However we can prove  some version of (\ref{er1})
 for sequences $f_n$, $cn^{-r}\le \|f_n\| \le Cn^{-r}$.

  For the problems of signal detection we say that maxisets $\gamma U$, $\gamma >0$, are perfect if the following two statements take place
\vskip 0.3cm
{\sl i.}  sequence of alternatives $f_n \in \mathbb{L}_2(0,1)$, $cn^{-r}\le \|f_n\| \le Cn^{-r}
$, is consistent iff  there are $\gamma U$, $\gamma > 0$, and sequence $f_{1n} \in \gamma U$, $c_2 n^{-r}\le \|f_{1n}\| \le C_2 n^{-r}
$, such that there holds
\begin{equation}
\label{ma1}
\| f_n \| =  \| f_{1n}\|  + \|f_n - f_{1n}\|.
\end{equation}
\vskip 0.3cm
{\sl ii.} sequence of alternatives $f_n \in \mathbb{L}_2(0,1)$, $cn^{-r}\le \|f_n\| \le Cn^{-r}
$, is inconsistent iff  for any $\gamma U$, for any sequence $f_{1n} \in \gamma U$, $c_2n^{-r}\le \|f_{1n}\| \le C_2n^{-r}
$, sequence of alternatives $f_n + f_{1n}$ is consistent and there holds
\begin{equation}
\label{ma2}
\|f_n + f_{1n}\| = \| f_n \| + \| f_{1n}\| + o(n^{-r})
\end{equation}
as $n \to \infty$.

For the problem of hypothesis testing, in {\sl ii.}, the sequence $f_{1n}$ should be such that $1 + f_n(s) + f_{1n}(s) > 0$ for all $s \in [0,1)$ and $\int_0^1 f_{1n}(s) \, ds =0$.

 As we know, a sequence of alternatives $f_{1n} \in \gamma U$, $c_2 n^{-r}\le \|f_{1n}\| \le C_2n^{-r}
$, $\gamma >0$, is consistent. Therefore {\sl i.} implies that, from any consistent  sequence of alternatives $f_n \in \mathbb{L}_2$, $cn^{-r}\le \|f_n\| \le Cn^{-r}
$, we can extract  a sufficiently smooth sequence of functions  $f_{1n} \in c_1U$, $c_2 n^{-r}\le \|f_{1n}\| \le C_2n^{-r}
$, responsible for consistency. Note that the sequence of alternatives $f_n$ is not necessary smooth and can be, for example, fast oscillating.
If we take a sequence inconsistent alternatives $f_n\in \mathbb{L}_2$, $c n^{-r}\le \|f_n\| \le C n^{-r}
$,(these functions could be unsmooth or fast oscillating), and add to functions $f_n$ any functions $f_{1n} \in c_1U$, $c_2 n^{-r} \le \|f_{1n}\| \le C_2n^{-r}
$,  (these functions can be considered as sufficiently smooth) then we get a consistent sequence of alternatives. Thus we can consider perfect maxisets as the kernels  generating all $n^{-r}$--consistent sequences.
 We can say that $n^{-r}$--consistent sequences of alternatives stringed in perfect maxisets.

 As we said the  maxisets are ambiguously defined. However the functions of any maxiset contains additive components from perfect maxisets.

For kernel--based test statistics and chi-squared test statistics, in the further  reasoning, we suppose  that the inconsistency of sequence of alternatives  $f_n$ in {\sl ii.} takes place  with arbitrary choice of windows width $h_n \asymp  n^{4r-2}$ and  with arbitrary choice of number of cells $k_n \asymp n^{2 - 4r}$ respectively. At the same time we suppose that the consistency takes place for some choice of windows width $h_n \asymp  n^{4r-2}$ and  for some choice of number of cells $k_n \asymp n^{2 - 4r}$ respectively.

 For any $\gamma >0$,  for any $ f \in \mathbb{L}_2(0,1)$ define the function $f_\gamma \in \gamma U$ such that $\|f_\gamma - f\| = \rho(f, \gamma U) = \inf \{ \|g - f\|, g \in \gamma U\}$.

The proof that test statistics $T_n$ satisfy {\sl i.}  in definition of perfect maxisets   is based on the choice of functions $f_{1n} = f_{\Pi_{i_n}}= \sum_{j=1}^{i_n} \theta_{jn} \phi_j$ for some sequence $i_n$. The functions $f_{\Pi_{i_n}}$ can be replaced with the function  $f_{n\gamma} = \sum_{j=1}^\infty\eta_{jn} \phi_j$. It is easy to show  that $\eta_{kn} = \theta_{kn}$ if $k^{2s}  \sum_{j=k}^\infty \theta^2_{jn} \le \gamma^2$. Hence $\eta_{kn} = \theta_{kn}$ for $k < \gamma^{1/s}  c^{-1/s} n^{r/s}$. This allows to hold the same reasoning for the sequence $f_{1n} = f_{n\gamma}$ as in the case $f_{1n} = f_{\Pi_{i_n}}$.

In particular all above mentioned test statistics satisfy as follows
\vskip 0.3cm
{\sl iii}. for any  $n^{-r}$-consistent sequence of alternatives $f_n$ there is $\gamma >0$ such that $f_{n\gamma}$ is $n^{-r}$-consistent.
\vskip 0.3cm
{\sl iv.} For any $\gamma > 0$, for any $n^{-r}$-inconsistent sequence of alternatives $f_n$ there holds $\|f_{n\gamma}\| = o(n^{-r})$.
\vskip 0.3cm
For  tests generated  $\mathbb{L}_2$- norms of kernel estimators and  tests generated quadratic forms of estimators of Fourier coefficients, the sequences $f_{n\gamma}$ have also the following property.

 We say that maxisets  $\gamma U$, $\gamma > 0$, are completely perfect for a sequence of test statistics $T_n$ if, for any $\epsilon > 0$ and any positive constants $c$  and $C$, $c < C$, there are  $\gamma_\epsilon$ and $n_\epsilon$ such that if sequence of alternatives $f_n \in \mathbb{L}_2$,   $c n^{-r}\le \|f_{n}\| \le Cn^{-r}
$, is consistent then, for any $n > n_\epsilon$, there hold
\begin{equation}
\label{uuu}
|\beta(K_n,f_n) - \beta(K_n,f_{n\gamma})| \le \epsilon
\end{equation}
and
\begin{equation}
\label{uu1}
|\beta(K_n,f_n-f_{n\gamma})| \ge 1 - \alpha - \epsilon,
\end{equation}
for all $\gamma > \gamma_\epsilon$.

Here $K_n$, $\alpha(K_n) = \alpha (1+ o(1))$ as $n \to \infty$, is a sequence of tests generated test statistics $T_n$.

By this definition we show that, for $n^{-r}$--consistent sequence of alternatives $f_n$, their projections $f_{n\gamma}$ on maxiset have almost the same type  II error probabilities and the influence of differences $f_n - f_{n\gamma}$ to the values of type II error probabilities is negligible.

For maxisets  of Cramer--von Mises tests we show that (\ref{uuu}) and (\ref{uu1}) hold if sequence $f_n$ satisfy some additional assumptions caused the requirement that $1 + f_{n\gamma}$ and $1 + f_n  - f_{n\gamma}$ should be densities.

\begin{remark} {\rm  In the framework of distance (semiparametric) approach asymptotic minimaxity of kernel-based tests, chi-squared tests  and test statistics generated quadratic forms of Fourier coefficients has been established for the wider sets of alternatives (see Ermakov \cite{er97, er03, er04, er11}). The asymptotic minimaxity of test statistics $T_n(\hat F_n)$ ($T_n(Y_n)$- for the problem of signal detection) has been proved for the sets of alternatives $\Psi_n = \{F : T_n(F) > b_n,\, F\,\, { \rm is\,\, c.d.f.}\}$ ($\Psi_n = \{f : T_n(f) > b_n, f \in L_2(0,1)\}$ respectively). The proof of results on maxisets can be treated as a search of the largest subset $V_n \subset \Psi_n$ such that alternatives $f_{i_n}$ in {\sl ii.} of definition of maxispace satisfies $T_n(f_{i_n}) > \tilde b_{i_n}$  (here we use the notation of problem of signal detection). Here $\tilde b_{i_n}$ are some constants defined by problem setup. Thus the role of  sets $V_n$ is somewhat blurred by the existence of   the larger sets $\Im_n$ of alternatives satisfying asymptotic minimaxity requirements. The notion of perfect maxisets and completely perfect maxisets emphasizes the role of maxisets. The perfect maxisets carry all information about $n^{-r}$--consistency and $n^{-r}$-inconsistency of sequences of alternatives}.\end{remark}
\section{Necessary conditions on consistency on the set $U$}
In all research on asymptotically nonparametric hypothesis testing with "small $\mathbb{L_2}$--ball removed" (see Ingster and Suslina \cite{ing02},    Yu. I. Ingster, T. Sapatinas, I. A. Suslina,  \cite{ing12}
  and
 Comminges and Dalalyan \cite{dal} (see also references therein)) the set $U$ is compact. Theorem \ref{qqqq} provided below shows that, if some assumption holds, this is necessary condition.

 We  consider the problem of signal detection in Gaussian white noise discussed in introduction. The problem will be explored in terms of sequence model.

 The stochastic differential equation (\ref{q1}) can be rewritten  in terms of  a sequence model for orthonormal system of functions $\phi_j$, $1 \le j < \infty$, in the following form
\begin{equation}\label{q2}
y_j = \theta_j + \frac{\sigma}{\sqrt{n}} \xi_j, \quad 1 \le j < \infty
\end{equation}
where $$y_j = \int \phi_j dY_n(t), \quad \xi_j = \int\,\phi_j\,dw(t) \mbox{ and}  \quad \theta_j = \int f\,\phi_j\,dt.$$ Denote $y = \{y_j\}_{j=1}^\infty$ and $\theta = \{\theta_j\}_{j=1}^\infty$.

In this notation the problem of hypothesis testing can be rewritten in the following form.
One needs to test the hypothesis $H_0 : \theta = 0$ versus alternatives $H_n : \theta \in \{ \theta : \|\theta\| \ge \rho_n, \theta \in U\}$ where $U$ is the ball with center  at zero  in functional space $\Im \subset l_2$.

We say that set $U$ is orthosymmetric if $\theta = \{\theta_j\}_{j=1}^\infty \in U$ and $|\eta_j| \le  |\theta_j|$ for all $j$ implies $\eta = \{\eta_j\}_{j=1}^\infty \in U$.
 \begin{thm} \label{qqqq} Suppose that set U is convex and orthosymmetric. Then there is consistent tests only if the set $U$ is compact.
 \end{thm}
Proof. If the set $U$ is convex and orthosymmetric, vectors $e_i$ in definition of maxisets have a simple form $e_i = \{e_{ij}\}_{j=1}^\infty$ with $e_{ij} = \delta_{ij}$ where  $\delta_{ij} =1$ if $i = j$ and $\delta_{ij} = 0$ if $i \ne j$.

If $U$ is not compact, this implies that there exists subsequence $i_l$   such that $d_{i_l} > c$ for all $l$. This implies that there is $n_0$ such that $e_{i_l} \in V_n$ for all $l > l(n_0)$.

Therefore, by Theorem 5.3 in Ermakov \cite{er15}, for any $n > n_0$ there does not exist uniformly consistent tests. This completes proof of Theorem \ref{qqqq}.
\section{Maxisets of quadratic test statistics}
We  consider the problem of signal detection in Gaussian white noise discussed in introduction. The problem will be explored in terms of sequence model.

 If $U$ is compact ellipsoid in Hilbert space, the asymptotically minimax test statistics are  quadratic forms
$$
T_n(Y_n) = \sum_{j=1}^\infty \kappa_{jn}^2 y_j^2 - \sigma^2 n^{-1} \sum_{j=1}^\infty \kappa_{jn}^2,
$$
with some specially defined coefficients $\kappa^2_{jn}$ (see Ermakov \cite{er90}).

If coefficients $\kappa_{jn}^2$ satisfy some regularity assumptions, the test statistics $T_n(Y_n)$ are asymptotically minimax for the  wider sets of alternatives
$$
{\rm H_n} : f \in Q_n(c)  = \{\, \theta:\, \theta= \{\theta_j\}_{j=1}^\infty,\,  A_n(\theta) > c \,\},
$$
with
$$
A_n(\theta) = n^{2}\,\sigma^{-4}\,\sum_{j=1}^\infty\, \kappa_{jn}^2\,\theta_j^2
$$
(see Ermakov \cite{er04}).

A sequence of tests $L_n, \alpha(L_n) = \alpha(1+ o(1))$, $0 <\alpha<1$, is called asymptotically minimax  if, for any sequence of tests $K_n, \alpha(K_n) \le \alpha,$ there holds
\begin{equation}\label{q4}
\liminf_{n\to \infty}(\beta(K_n,Q_n(c)) - \beta(L_n,Q_n(c))) \ge 0.
\end{equation}
Sequence of test statistics $T_n$ is asymptotically minimax if the  tests generated test statistics $T_n$ are asymptotically minimax.

Section goal is to point out maxisets for test statistics $T_n(Y_n)$ with coefficients $\kappa^2_{jn}$ satisfying some regularity assumptions.

Assume that
the coefficients $\kappa_{jn}^2, 1\le j < \infty,$ satisfy the  following assumptions.

\noindent{\bf A1.} For each $n$ the sequence $\kappa^2_{jn}$ is decreasing.

\noindent{\bf A2.} There are positive constants $C_1,C_2$ such that, for each $n$, there holds
\begin{equation}\label{q5}
 C_1 < A_n = \sigma^{-4}\,n^2\,\sum_{j=1}^\infty \kappa_{jn}^4 < C_2.
 \end{equation}
  Denote $$k_n = \sup\Bigl\{k: \sum_{j < k} \kappa^2_{jn} \le \frac{1}{2} \sum_{j =1}^\infty \kappa^2_{jn} \Bigr\}. $$
\
\noindent{\bf A3.} For any $\delta$,  $0 <\delta < 1/2,$ there holds
\begin{equation}\label{q7}
\lim_{n\to\infty} \sup_{\delta k_n < j < \delta^{-1}k_n} \left|\frac{\kappa_{j+1,n}^2}{\kappa_{j,n}^2} -1\right| = 0.
\end{equation}
\noindent{\bf A4.} For any $\delta > 0$ and any $\delta_1$, $0 \le \delta_1 < 1$, there are $C_1$  and $C_2<1$ such that
 \begin{equation}\label{q7s}
C_1 <  \frac{\kappa^2_{(1+\delta)k_n,n}}{\kappa^2_{(1-\delta_1)k_n,n}} < C_2.
\end{equation}
\noindent{\bf A5.}
\begin{equation}\label{q8}
\lim_{\delta\to 0} \lim_{n\to \infty}\frac{\sum_{\delta k_n < j < \delta^{-1}k_n} \kappa_{jn}^2}{\sum_{j=1}^\infty \kappa_{jn}^2}  =1
\end{equation}
and
\begin{equation}\label{q9}
\lim_{\delta\to 0} \lim_{n\to \infty}  A_n^{-1} n^2\sum_{\delta k_n < j < \delta^{-1}k_n} \kappa_{jn}^4 =1
\end{equation}

{\sl Example}.  Let
$$
\kappa^2_{jn} = n^{-1/(2\gamma)}\frac{n^{-1} j^{-\gamma}}{j^{-\gamma}  + n^{-1}}, \quad \gamma > 0.
$$
Then A1 -- A5 hold.

Denote $s = \frac{r}{2 -4r}$. Then $r = \frac{2s}{1 + 4s}$.
\begin{thm}\label{tq1} Assume {\rm A1-A5}. Then the balls  $\mathbb{\bar B}^s_{2\infty}(P_0)$ are maxisets for the test statistics $T_n(Y_n)$ with $k_n \asymp n^{2-4r}=n^{\frac{2}{1+4s}}$. \end{thm}
\begin{thm}\label{tq3} Assume {\rm A1-A5}. Then the balls  $\mathbb{\bar B}^s_{2\infty}(P_0)$ are perfect maxisets.
\end{thm}
The balls $\mathbb{\bar B}^s_{2\infty}(P_0)$ in Theorems \ref{tq1} and  \ref{tq3} can be replaced with any  ball in $\mathbb{\bar B}^s_{2\infty}$ generated equivalent norm.

\begin{thm}  \label{tq4}  Assume {\rm A1-A5}. Then the balls  $\mathbb{\bar B}^s_{2\infty}(P_0)$ are completely perfect maxisets.
\end{thm}
\begin{remark}\label{rem1}{\rm Let $\kappa^2_{jn} = 0$ for $j > l_n$ and  let $\kappa_{jn}^2 > 0$ for $j \le l_n$ with $l_n \to \infty$ as $n \to \infty$. The analysis of the proofs of Theorems   \ref{tq1} and \ref{tq3} shows that Theorems \ref{tq1} - \ref{tq4}  remain valid for this setup  if we make the following changes in A1 -- A4.
We put $k_n = l_n$. We replace $\delta^{-1} k_n$ with $(1-\delta)  k_n$ in (\ref{q7}), (\ref{q8}), (\ref{q9}) and replace (\ref{q7s}) with
 \begin{equation}\label{q7v}
C_2(\delta) > \frac{\kappa^2_{(1-\delta)k_n,n}}{\kappa^2_{n}} > C_1(\delta) > 0
\end{equation}
with $\kappa^2_n = \kappa^2_{k_n/2,n}$. Here $0 <  \delta < 1$. In the corresponding version of Theorem \ref{tq1} one needs also to require additionally $k_n > c_0\,n^{2-4r}$ for some $c_0=c_0(P_0) > 0$. The differences in the reasoning are the same as in the proofs of Theorems \ref{tk1} - \ref{tk4} of the next section}.\end{remark}
\section{Maxisets of kernel-based tests}
We  explore the problem of signal detection of previous section and suppose additionally that function $f$ belongs to $\mathbb{L}_2^{per}(\mathbb{R}^1)$ the set of 1-periodic functions such that $f(t) \in \mathbb{L}_2(0,1), t \in [0,1)$. This allows to extend our model on real line $\mathbb{R}^1$ putting $w(t+j) = w(t)$ for all integer $j$ and $t \in [0,1)$ and  to write the forthcoming integrals over all real line.

Define kernel estimator
\begin{equation}\label{yy}
\hat{f}_n(t) = \frac{1}{h_n} \int_{-\infty}^{\infty} K\Bigl(\frac{t-u}{h_n}\Bigr) dY_n(u), \quad t \in (0,1),
\end{equation}
where $h_n$ is a sequence of positive numbers, $h_n \to 0$ as $n \to 0$. The kernel $K$ is bounded function such that the support of $K$ is contained in $[-1,1]$, $K(t) = K(-t)$ for $t \in R^1$ and $\int_{-\infty}^\infty K(t) dt = 1$.

In (\ref{yy}) we suppose that, for any $v, 0 <v  < 1$, we have
$$
\frac{1}{h_n} \int_{1}^{1+v} K\Bigl(\frac{t-u}{h_n}\Bigr)\,dY_n(u) = \frac{1}{h_n} \int_{0}^{1+v} K\Bigl(\frac{t-1-u}{h_n}\Bigr)\,f(u)\,du
$$
$$
 + \frac{\sigma}{\sqrt{n}h_n} \int_{0}^{1+v} K\Bigl(\frac{t-1-u}{h_n}\Bigr)\, dw(u)
$$
and
$$
\frac{1}{h_n} \int_{-v}^{0} K\Bigl(\frac{t-u}{h_n}\Bigr)\,dY_n(u) = \frac{1}{h_n} \int_{1-v}^{1} K\Bigl(\frac{t-u+1}{h_n}\Bigr)\,f(u)\,du
 $$
 $$
 + \frac{\sigma}{\sqrt{n} h_n} \int_{1-v}^{1} K\Bigl(\frac{t-u+1}{h_n}\Bigr)\, dw(u).
$$
For hypothesis testing we implement the kernel-based tests (see Bickel and Rosenblatt \cite{bic}) with the test statistics $$T_n(Y_n) =nh_n^{1/2}\sigma^{-2}\kappa^{-1} (\|\hat f_{h_n}\|^2- \sigma^2(nh_n)^{-1}\|K\|^2)$$  where
$$
\kappa^2 = 2 \int \Bigl(\int K(t-s)K(s) ds\Bigr)^2\,dt.
$$
\begin{thm}\label{tk1}
  Balls  $\mathbb{B}^{s}_{2\infty}(P_0)$ in Besov space $\mathbb{B}^{s}_{2\infty}$ with $s =\frac{r}{2-4r}$ are maxisets for kernel-based tests with $h_n \asymp n^{4r-2}=n^{\frac{-2}{1+4s}}$ and $h_n < c_0\,n^{4r-2}$ for some $c_0 = c_0(P_0) >0$.
\end{thm}
\begin{thm}\label{tk3} Balls  $\mathbb{B}^{s}_{2\infty}(P_0)$ in Besov space $\mathbb{B}^{s}_{2\infty}$ with $s =\frac{r}{2-4r}$are perfect maxisets.
\end{thm}
\begin{thm}  \label{tk4}  Balls  $\mathbb{B}^{s}_{2\infty}(P_0)$ in Besov space $\mathbb{B}^{s}_{2\infty}$ with $s =\frac{r}{2-4r}$ are completely perfect maxisets.
\end{thm}
\section{Maxisets of $\chi^2$-tests}
Let $X_1,\ldots,X_n$ be i.i.d.r.v.'s having c.d.f. $F(x)$, $x \in (0,1)$.
Let c.d.f. $F(x)$ has a density $1 + f(x) = dF(x)/dx, x \in (0,1), f \in L_2^{per}(0,1)$.
We explore the problem of testing hypothesis (\ref{i1}) and (\ref{i2b}) discussed in introduction.

Let $\hat F_n(x)$ be empirical c.d.f. of $X_1,\ldots,X_n$.

Denote  $\hat p_{in} = \hat F_n((i+1)/k_n) - \hat F_n(i/k_n), 1 \le i \le k_n$.

The test statistics of $\chi^2$-tests equal
$$
T_n(\hat F_n) =  k_n n \sum_{i=1}^{k_n} (\hat p_{in}  - 1/k_n)^2.
$$
\begin{thm}\label{chi1}
  Balls  $\mathbb{B}^{s}_{2\infty}(P_0)$  in Besov spaces $\mathbb{B}^{s}_{2\infty}$ with $s =\frac{r}{2-4r}$ are maxisets for $\chi^2$-tests with the number of cells $k_n \asymp n^{2- 4r}=n^{\frac{2}{1+4s}}$ and $k_n > c_0\,n^{2- 4r}$ for some $c_0= c_0(P_0) >0$.
\end{thm}
\begin{thm}\label{chi3} The  balls  $\mathbb{B}^{s}_{2\infty}(P_0)$ in Besov spaces $\mathbb{B}^{s}_{2\infty}$ with $s =\frac{r}{2-4r}$ are  perfect maxisets.
\end{thm}
{\sl Discussion}
The definition of $\chi^2$ - tests is based on indicator functions. Thus $\chi^2$ - tests should detect well distribution functions with stepwise densities.
Besov spaces $\mathbb{B}^{s}_{2\infty}, s \ge 1,$  do not contain stepwise  functions. It seems strange.

Let us consider $\chi^2$ - test with $k_n = 2^{l_n}, l_n \to \infty$ as $n \to \infty$. Then  $\chi^2$ - test statistics admit representation
$$
T_n(\hat F_n) =  k_n n \sum_{i=1}^{l_n}\sum_{j=1}^{2^i} \hat \beta_{ij}^2,
$$
with
$$
\hat \beta_{ij} = \frac{1}{n}\sum_{m=1}^n \phi_{ij}(X_m),
$$
where $\phi_{ij}$ are functions of Haar orthogonal system, $\phi_{ij}(x) =  2^{i/2} \phi(2^ix-j)$ with $\phi(x) = 1$ if $x \in (0,1/2)$, $\phi(x) =-1$  if $x \in (1/2,1)$ and $\phi(x)=0$ otherwise.

 Implementing the same reasoning as in the case quadratic test statistics and using Theorem \ref{chi2} given below, we get that $\chi^2$ - test statistics have maxisets
$$
\bar B^s_{2\infty}(P_0) = \Bigl\{f : f = 1+ \sum_{k=1}^\infty\sum_{j=1}^{2^k}\beta_{kj}\phi_{kj},\,\,\,  \sup_{\lambda>0} 2^{2\lambda s} \sum_{k>\lambda}^\infty\sum_{j=1}^{2^k}\beta^2_{kj} \le P_0\Bigr\}.
$$
This statement is true as well.

Suppose  function $f$ is sufficiently smooth and $\beta_{kj}$ are Fourier coefficients of $f$ for Haar orthogonal system. Since $\beta_{kj} = 2^{-k/2}\frac{df}{dx}(j2^{-k})(1 +o(1))$ as $k \to \infty$, then
$$
\sum_{j=1}^{2^k}\beta^2_{kj} = C2^{-k/2}\int \Bigl(\frac{df}{dx}\Bigr)^2\, dx (1+ o(1)).
$$
Thus we see that $f$ does not belong to $\mathbb{B}^s_{2\infty}, s >1,$ for such a setup.

Kernel-based tests also detect stepwise densities well. However these densities  does not also belong the maxispaces of kernel-based tests.
\section{Maxisets of Cramer -- von Mises tests}
We  consider Cramer -- von Mises test statistics as functionals
$$
T^2(\hat{F}_n -F_0) = \int_0^1 (\hat{F}_n(x) - F_0(x))^2\, d F_0(x)
$$
depending on empirical distribution function $\hat F_n$.
Here $F_0(x)=x, x \in (0,1)$.

The functional $T$ is the norm on the set of differences of distribution functions. Therefore we have
\begin{equation}\label{cr1}
T(\hat{F}_n - F_0) - T(F - F_0) \le T(\hat{F}_n - F) \le T(\hat{F}_n - F_0)  + T(F- F_0).
\end{equation}
Hence  it is easy to see that sequence of  alternatives $F_n$ is consistent iff
\begin{equation}\label{cru}
nT^2(F_n - F_0) > c  \quad \hbox {\rm for all}\quad n > n_0
\end{equation}
This allows to search for the maxiset as the largest convex set $U \subset \mathbb{L}_2(0,1)$ satisfying the following conditions
\vskip 0.3cm
\noindent{\sl i.} for all $f = \frac{d(F-F_0)}{dx} \in U$ such that $cn^{-r} < \|f\| < C n^{-r}$, there holds \begin{equation}\label{cr2}
\sqrt{n} T(F - F_0) > c
\end{equation}
{\sl ii.} for any $f \notin \lambda U$ for all $\lambda >0$, $\int_0^1 f(x)dx = 0$, there are sequences $i_n,j_n$ such that $cj_n^{-r}\le \|\tilde f_{i_n}\| \le C j_n^{-r}$ and, if $1 + f_n(x) > 0$ for all $x \in (0,1)$, then there holds
\begin{equation}\label{}
\lim_{n \to \infty} j_n^{1/2} T(\tilde F_{i_n} - F_0) = \infty
\end{equation}
with $\frac{d\tilde F_{i_n}}{dx} - 1 = \tilde f_{i_n}$.
\begin{thm}\label{tom1}
The balls  $\mathbb{B}^s_{2\infty}(P_0)$ with $s = \frac{2r}{1 - 2r}$,  $r = \frac{s}{2+2s}$, are maxisets for Cramer -- von Mises test statistics. Here  the orthonormal functions $\phi_j(x) = \sqrt{2}\cos(\pi j x)$, $x \in [0,1]$, $1 \le j < \infty$.
\end{thm}
Here the balls $\mathbb{B}^s_{2\infty}(P_0)$ are defined (\ref{vv}).
\begin{thm}\label{tom2} The balls  $\mathbb{B}^s_{2\infty}(P_0)$ with $s = \frac{2r}{1 - 2r}$ are perfect maxisets for Cramer -- von Mises test statistics. Here  the orthonormal functions $\phi_j(x) = \sqrt{2}\cos(\pi j x)$, $x \in [0,1]$, $1 \le j < \infty$.
\end{thm}
 For Cramer- von Mises tests we made additional assumptions in definition of completely perfect maxisets. We fix any $\delta > 0$. We state  that (\ref{uuu}) and (\ref{uu1}) holds for the sequences $f_n$ such that  B1  is fulfilled.
\vskip 0.25cm
{\bf B1}. For all $x \in (0,1)$ and all $\gamma > \gamma_0$ there hold $1 + f_n(x) > \delta$,  $1+ f_{n\gamma}(x) > \delta$ and $1+ f_n(x)-f_{n\gamma}(x) > \delta$.
\begin{thm}\label{tom2b} The balls $\mathbb{B}^s_{2\infty}(P_0)$ with $s = \frac{2r}{1 - 2r}$ are completely perfect maxisets for Cramer -- von Mises test statistics.
\end{thm}
\section{  Asymptotically minimax tests for maxisets}
Let we observe a random process $Y_n(t)$, $t \in [0,1)$, defined by the stochastic differential equation (\ref{q1}) with unknown signal $f$.

Our goal is to  point out asymptotically minimax tests for the problem testing of the hypothesis
${\rm H_0\,:} \quad f(t) = 0$,\,\,\, $t \in [0,1) $,
versus the alternatives
$$ {\rm H_n\,:}\quad   \|f\|^2 > \rho_n \asymp n^{-\frac{4s}{1+4s}}
$$
if a priori information is provided that
$
f \in \mathbb{\bar B}^s_{2\infty}(P_0).
$

Denote $V_n = \{f:\, \|f\|^2 \ge \rho_n,\, f \in  \mathbb{\bar B}_{2\infty}^s(P_0)\}$.

Note that, for Besov balls
$$
\mathbb{\tilde B}^s_{2\infty}(P_0) =  \Bigl\{f : f =  \sum_{k=1}^\infty\sum_{j=1}^{2^k}\theta_{kj}\phi_{kj},\,\,\,  \sup_{k} 2^{2k s} \sum_{j=1}^{2^k}\theta^2_{kj} \le P_0\Bigr\}
$$
provided in terms of wavelet functions, asymptotically minimax tests have been established Ingster and Suslina \cite{ing02}. Here the assignment of Besov ball is different.

 In estimation, for  Besov balls $\mathbb{\bar B}^s_{2\infty}(P_0)$ we get that penalized maximum likelihood estimators are asymptotically minimax \cite{er17}. This illustrates the role of such a priori information in statistical inference.

The proof, in main features, repeats the reasoning in Ermakov \cite{er90}. The main difference   is the solution of another extremal problem
caused by another definition of sets of alternatives. Other differences have technical character and are also caused the differences of definition of sets of alternatives.

The results will be provided in terms of sequence model (see section 3).

Define $k=k_n$ and $\kappa^2 = \kappa_n^2$ as the solution of two equations
\begin{equation}\label{i2}
\frac{1}{2s}k_n^{1+2s}\kappa^2_n = P_0
\end{equation}
and
\begin{equation}\label{i3}
k_n\kappa_n^2 + k_n^{-2s} P_0  = \rho_n.
\end{equation}
Denote
$
\kappa_j^2 =  \kappa_n^2$, for $1 \le j \le k_n
$ and $\kappa_j^2= 2s P_0 j^{-2s-1}$, for $ j > k_n.$

Define test statistics
$$
T^a_n(Y_n) =   \sigma^{-2}n \sum_{j=1}^\infty \kappa_j^2 y_j^2.
$$
and put
$$
A_n = \sigma^{-4}n^{2}\sum_{j=1}^\infty \kappa_j^4,
$$
$$
C_n =\sigma^{-2} n\rho_n.
$$
For type I error probabilities $\alpha, 0 < \alpha< 1$, define the critical regions
$$
S^a_n == \{ y: \, (T^a_n(y) - C_n)(2A_n)^{-1/2} > x_\alpha\}
$$
with $x_\alpha$ defined by equation $\alpha = 1 - \Phi(x_\alpha).$
\begin{thm}\label{t1} Let
\begin{equation}\label{i4}
0 <\liminf_{n \to \infty} A_n \le \limsup_{n\to \infty} A_n< \infty.
\end{equation}
Then the tests $L^a_n$ with critical regions $S^a_n$ are asymptotically minimax with $\alpha(L^a_n) = \alpha(1+o(1))$ and
 \begin{equation}\label{i5}
\beta(L_n^a,V_n) = \Phi(x_\alpha - (A_n/2)^{1/2})(1 + o(1))
\end{equation}
as $n \to \infty$.
\end{thm}
\noindent{\sl Example.} Let $\rho_n = R (\sigma^2/n)^{\frac{4s}{1+4s}}(1+o(1))$ as $n \to \infty$. Then
$$
A_n = \sigma^{-4}n^2\rho_n^{\frac{1+4s}{2s}}\frac{8s^2}{(1+4s)(1+2s)}((1+2s)P_0)^{-1/2s}(1 + o(1))
$$
$$= R^2\frac{8s^2}{(1+4s)(1+2s)}((1+2s)P_0)^{-1/2s}(1 + o(1)).
$$
Ingster, Sapatinas, Suslina \cite{ing12} and Laurent, Loubes, Marteau \cite{la} have explored the problem of signal detection for linear inverse ill-posed problems. The setup was treated in terms of sequence model
$$
y_j = \lambda_j\theta_j + \frac{\sigma}{\sqrt{n}} \xi_j, \quad 1 \le j < \infty
$$
where $\xi_j$ are i.i.d.r.v.'s having standard normal distribution and $\lambda_j$ is sequence of eigenvalues of linear operator.

It is easy to see that,  if $|\lambda_j| \asymp j^{-\gamma}$, then the maxisets for tests statistics defined as quadratic forms of $y_j$, $1 \le j < \infty$, are the balls in $\mathbb{B}^s_{2\infty}$  with $r = \frac{2s}{1+ 4s +4\gamma}$.
Thus it is of interest to point out asymptotically minimax test statistics for the problem of testing of hypothesis ${\rm H_0}\, :\, \theta = 0$ versus alternatives ${\rm H_n} \,:\, \theta \in V_n$.

Define test statistics
$$
T^a_n(Y_n) =   \sigma^{-2}n \sum_{j=1}^\infty \kappa_j^2 y_j^2,
$$
with $\kappa_j^2$ defined the equations
$\kappa_j^2 = a\lambda_j^{-2}$ for $j \le k_n$ and $\kappa_j^2 = 2s P_0 \lambda_j^2 j^{-1-2s}$ for $j > k_n$, where constants $a= a_n$ and $k_n$ are the solutions of  equations
$$
a_n \sum_{j=1}^{k_n} \lambda_j^{-4} + P_0 k_n^{-2s} = \rho_n(1 + o(1)) \quad \mbox{\rm and} \quad a_n\lambda_{k_n}^{-4} = 2s P_0 k_n^{-1-4s}(1+o(1)).
$$
In this notation  the definition of $A_n$ and the critical regions $S^a_n$ is the same as in Theorem \ref{t1}.
\begin{thm}\label{t2} Let  $|\lambda_j| \asymp j^{-\gamma}$. Then for the above setup and for above notation the statement of Theorem \ref{t1} holds.\end{thm}
\noindent{\sl Example.} Let $\lambda_j^2 = A j^{-2\gamma}$ and let $\rho_n \asymp n^{\frac{-4s}{1+4s+4\gamma}}$. Then
$$
A_n = \sigma^{-4} n^2 \rho_n^{\frac{1+4s+4\gamma}{2s}} A^2
 \frac{8 s^2 (1+ 4\gamma)}{(1 +2s + 4\gamma)(1 + 4s + 4\gamma)}\Bigl(\frac{1+2s +4\gamma}{1+4\gamma }P_0\Bigr)^{-\frac{1 + 4\gamma}{2s}}(1+o(1)).
 $$
 Proof of Theorem \ref{t2} is akin to that of Theorem \ref{t1} and is omitted.
 \section{Proof of Theorems}
 \subsection{Proof of Theorems of section 4}
 {\sl Proof of Theorem  \ref{tq1}. Sufficiency}. The proof is based on the inequality (\ref{u2}) defining the rate of consistency  and on the relation (\ref{aq6}) that balances   the contribution of bias and stochastic part of test statistics $T_n(Y_n)$.  This two relations assign in Theorem \ref{tq1} two parameters:   the limitation $k_n \asymp n^{2-4r}$ on coefficients $\kappa_{jn}^2$ and the order of decreasing of the tail $\theta=\{\theta_j\}_{j=1}^\infty \in \mathbb{\bar B}^s_{2\infty}$.

The reasoning is based on Theorem \ref{tq2} on asymptotic minimaxity of test statistics $T_n$.
\begin{thm}\label{tq2}
Assume {\rm A1-A5}. Then sequence of tests $K_n(Y_n) = \chi\{n^{-1}T_n(Y_n) > (2A_n)^{1/2} x_\alpha\}$ is asymptotically minimax for the sets of alternatives $Q_n(c)$.

There holds
\begin{equation}\label{aq1}
\beta(K_n,\theta) = \Phi(x_{\alpha} - A_n(\theta)(2A_n)^{-1/2})(1+o(1))
\end{equation}
uniformly in all $\theta$ such that $A_n(\theta)< C$. Here $x_\alpha$ is defined by the equation $\alpha = 1 - \Phi(x_\alpha)$.
\end{thm}
A version of Theorem \ref{tq2} for the model
$$
dY(t) = f(t)\,dt + \frac{\sigma}{\sqrt{n}} h(t)\,dw(t), \quad t \in (0,1),
$$
with  heteroscedastic white noise $ h \in L_2 (0,1)$ has been proved in Ermakov \cite{er03}.

{\sl Proof of Theorem \ref{tq2}}. Theorem \ref{tq2} and its version for Remark \ref{rem1} setup can be deduceded straightforwardly from Theorem 1 in Ermakov \cite{er90}. The lower bound follows from  Theorem 1 in \cite{er90}.

The upper bound follows from the following reasoning.  We have
\begin{equation}\label{aq2}
\begin{split}&
 \sum_{j=1}^\infty \kappa_{jn}^2 y_j^2  = \sum_{j=1}^\infty \kappa_{jn}^2 \theta_{jn}^2 + 2\frac{\sigma}{\sqrt{n}}\sum_{j=1}^\infty \kappa_{jn}^2 \theta_{jn} \xi_j +
\frac{\sigma^2}{n} \sum_{j=1}^\infty \kappa_{jn}^2 \xi_j^2\\&
= J_{1n} + J_{2n} + J_{3n},
\end{split}
\end{equation}
with
\begin{equation}\label{aq3}
\mathbf{E} [J_{3n}] = \frac{\sigma^2}{n} \rho_n,
\quad
\mathbf{Var} [J_{3n}] = 2 \frac{\sigma^4}{n^4} A_n,
\end{equation}
\begin{equation}\label{aq5}
\mathbf{Var} [J_{2n}] = \frac{\sigma^2}{n}  \sum_{j=1}^\infty \kappa_{jn}^4 \theta_{jn}^2 \le \frac{\sigma^2\kappa^2}{n}  \sum_{j=1}^\infty \kappa_{jn}^2 \theta_{jn}^2.
\end{equation}
By Chebyshov inequality, it follows from (\ref{aq2}) - (\ref{aq5}), that, if $A_n n^{-2} = o(J_{1n})= o\Bigl(\sum_{j=1}^\infty \kappa_{jn}^2 \theta_{jn}^2\Bigr)$ as $ n \to \infty$, then $\beta(L_n,\theta_n) \to 0$ as $n \to \infty$. Thus it suffices to explore the case
\begin{equation}\label{aq6}
A_n^2 \asymp n^2\sum_{j=1}^\infty \kappa_{jn}^2 \theta_{jn}^2.
\end{equation}
If (\ref{aq6}) holds, then implementing the reasoning of Theorem 1 in \cite{er90} we get that (\ref{aq1}) holds. This completes the proof of Theorem \ref{tq2}.

Let $\theta= \{\theta_j\}_{j=1}^\infty \in \mathbb{\bar B}^s_{2\infty}$.

Denote $\kappa^2 = \kappa^2_{k_nn}$.
Note that A1, A2 and A4 imply that
\begin{equation}\label{u1}
\kappa^4 \asymp n^{-2}k_n^{-1}.
\end{equation}
Without loss of generality, we can suppose that
 $||\theta||^2  \asymp n^{-2r}$.

 Then there is $k_n = Cn^{2-4r}(1 + o(1))$  such that
\begin{equation}\label{u2}
k_n^{2s}\sum_{j=1}^{k_n} \theta_j^2(1 + o(1))  = C_1 n^{2r}\sum_{j=1}^{k_n} \theta_j^2 > C_0
\end{equation}
where constants  $C$, $C_0$,  $C_1$  do not depend on $n$.

Otherwise, there is $C_3$ such that, for any $C_2$ and $k_n = C_2 n^{2-4r}(1 + o(1))$, we get
\begin{equation}\label{u3}
n^{2r}\sum_{j=k_n}^\infty \theta_j^2 > C_3
\end{equation}
that implies $\theta \notin \mathbb{\bar B}^s_{2\infty}$.

By  $||\theta||^2  \asymp n^{-2r}$ and (\ref{u1}), (\ref{u2}) together,  we get
\begin{equation}\label{u4}
n^2 \sum_{ j=1}^\infty \kappa_j^2\theta_j^2 \asymp n^2\kappa^2 \sum_{j=1}^\infty\theta_j^2 \asymp n^{1-2r} k_n^{-1/2} \asymp 1.
\end{equation}
It remains  to implement asymptotically minimax  Theorem \ref{tq2}.

{\sl Proof of necessary condition}. Suppose the opposite.
Then there are $\theta = \{\theta_j\}_{j=1}^\infty$, $\theta \notin \Im$,  and a sequence $m_l, m_l \to \infty$ as $l \to \infty$, such that
\begin{equation}\label{u5}
m_l^{2s} \sum_{j=m_l}^\infty \theta_j^2 = C_l
\end{equation}
with $C_l \to \infty$ as $l \to \infty$.

It is clear that  we can define a sequence $m_l$ such that
\begin{equation}\label{u6}
m_l^{2s} \sum_{j=m_l}^{2m_l} \theta_j^2 > \delta C_l, \quad 0 < \delta < 1/2,
\end{equation}
 where $\delta$ does not depend on $l$.

 Otherwise, we have
 $$
 2^{2s(i-1)} m_l^{2s} \sum_{j=2^{i-1}m_l}^{2^{i}m_l}  \theta_j^2  < \delta  C_l
 $$
 for  all $i = 1,2,\ldots$,  that implies that the left hand-side of
 (\ref{u5}) does not exceed $2\delta C_l$.

Define a sequence $\eta_l = \{\eta_{jl}\}_{j=1}^l$ such that
$\eta_{jl} = 0$ if $j  < m_l$ and $\eta_{jl} = \theta_j$ if $j \ge m_l$.

For alternatives $\eta_l$ we define sequence $n= n_l$ such that
\begin{equation}\label{u7}
n_l \asymp C_l^{-1/(2r)} m_l^{s/r} = C_l^{-1/(2r)} m_l^{\frac{1}{2 - 4r}}.
\end{equation}
Then
\begin{equation}\label{u5b}
\|\eta_l\|^2 \asymp m_l^{-2s} C_l \asymp n_l^{-2r}.
\end{equation}
Since sequence $\kappa^2_{jn_l}$ is decreasing and (\ref{u6}) holds, by (\ref{q7s}), we have
\begin{equation}\label{u5g}
\sum_{j=1}^\infty  \kappa^2_{jn_l} \eta^2_{jl} \asymp \kappa^2_{n_l} \sum_{j=m_l}^{2m_l} \eta^2_{jn_l}.
\end{equation}
Therefore $k_{n_l}  \asymp  m_l$. Denote $k_l = 2m_l$.

Then
\begin{equation}\label{u8}
k_l^{2s} \sum_{j=k_l/2}^{k_l} \eta_{jl}^2 \asymp C_l.
\end{equation}
Hence
\begin{equation}\label{u9}
k_l^{2s} n_l^{-2r} =  k_l^{\frac{2r}{2 - 4r}} n_l^{-2r}\asymp C_l.
\end{equation}
Therefore we get
\begin{equation}\label{u10}
k_l^{1/2} \asymp C_l^{(1-2r)/2}n_l^{1-2r}.
\end{equation}
By (\ref{u5}), (\ref{u6}) and A3, we get
\begin{equation}\label{u11}
\sum_{j=k_l/2}^{k_l} \kappa_{jn_l}^2\eta_{jl}^2 \asymp \sum_{j=1}^{\infty} \kappa_{jn_l}^2\eta_{jl}^2.
\end{equation}
Using (\ref{u1})  and (\ref{u10}), we get
\begin{equation}\label{u12}
n_l^2 \sum_{j=k_l/2}^{k_l} \kappa_{jn_l}^2\eta_{jl}^2 \asymp n_l k_l^{-1/2}\sum_{j=1}^{k_l} \eta_{jl}^2 \asymp n_l^{1-2r}k_l^{-1/2} \asymp C_l^{-(1-2r)/2}.
\end{equation}
By (\ref{u5g}) and Theorem \ref{tq2}, (\ref{u11}) and (\ref{u12}) imply inconsistency of   sequence of alternatives $\eta_l$.

{\sl Proof of Theorem \ref{tq3}}.  The reasoning  is based on Lemmas \ref{ld1} -- \ref{ld7}. Statement {\sl i.} follows from Lemmas \ref{ld4} and \ref{ld6}.  Statement {\sl ii.} follows from Lemmas \ref{ld5} and \ref{ld7}.

\begin{lemma} \label{ld1} Let $cn^{-r}\le \|f_n\| \le Cn^{-r}
$ and $f_n \in c_1 U$. Then, for $k_n = C_1 n^{2 - 4r}(1  +o(1)) = C_1 n^{\frac{2}{1 + 4s}}(1+o(1)) $ with $C_1 > c/(2c_1)$, there holds
\begin{equation}\label{d1}
\sum_{j=1}^{k_n} \theta_{jn}^2 > \frac{c}{2} n^{-2r}.
\end{equation}
\end{lemma}
Proof. If $k_n^{2s} = C_1^{2s} n^{2r}(1 + o(1))$ and $f_n \in c_1 U$, then we have
\begin{equation}\label{d2}
k_n^{2s} \sum_{j=k_n}^\infty \theta_{jn}^2 = C_1^{2s} n^{2r}  \sum_{j=k_n}^\infty \theta_{jn}^2(1 + o(1))\le c_1.
\end{equation}
Hence
\begin{equation}\label{d3}
\sum_{j=k_n}^\infty \theta_{jn}^2 \le c_1 C_1^{-2s} n^{-2r} \le \frac{c}{2} n^{-2r}.
\end{equation}
Therefore (\ref{d1}) holds.
\begin{lemma}\label{ld2} Let sequence $f_n$ be $n^{-r}$ -inconsistent for $T_n$ with $k_n \asymp n^{2 - 4r}$. Then, for any $c$, there holds
\begin{equation}\label{d4}
k_n^{2s} \sum_{j=1}^{ck_n} \theta_{jn}^2 \asymp n^{2r}  \sum_{j=1}^{ck_n} \theta_{jn}^2 = o(1).
\end{equation}
\end{lemma}
Here the summation is over  all $1 \le j < ck_n$.  In what follows, we shall use this notation as well.

Proof. Suppose opposite. Then, by A4 and (\ref{u1}), we have
\begin{equation}\label{d5}
n^{2}  \sum_{j=1}^{ck_n} \kappa^2_{jn}\theta_{jn}^2 \asymp n^2  \kappa^2 \sum_{j=1}^{ck_n} \theta_{jn}^2 \asymp n^{2r} \sum_{j=1}^{ck_n} \theta_{jn}^2
\end{equation}
By Theorem \ref{tq2}, this implies (\ref{d4}).

\begin{lemma}\label{ld3} For any $c$ and $C$ there is $\gamma$ such that if $\| f_n \| \le C n^{-r}$ and
$
f_n = \sum_{j =1}^{ck_n} \theta_{jn} \phi_j
$
then  $f_n \in\gamma U$.
\end{lemma}
Proof. We have
\begin{equation}\label{d6}
k_n^{2s} \sum_{j=1}^{ck_n} \theta_{jn}^2 \le  C_1 n^{2r}  \sum_{j=1}^{\infty} \theta_{jn}^2 < C.
\end{equation}
This implies Lemma \ref{ld3}.
\begin{lemma}\label{ld4} Let  (\ref{ma1}) hold. Then sequence $f_n$ is $n^{-r}$-consistent.
\end{lemma}
Let $f_n = \sum_{j=1}^\infty \theta_{jn}  \phi_j$ and let  $$f_{1n} = \sum_{j=1}^\infty \eta_{jn}  \phi_j, \quad f_n - f_{1n} =\sum_{j=1}^\infty \zeta_{jn}  \phi_j.$$
For any $\delta > 0$ there is $c$ such that
\begin{equation}\label{uh11}
\sum_{j>ck_n} \eta_{jn}^2 < \delta n^{-2r}
\end{equation}
for each $f_{1n} \in c_1U$, $\| f_n \| \le C_2 n^{-r}$ .

We have
\begin{equation}\label{dub2}
\begin{split}&
J_n=\left| \sum_{j > ck_n} \theta_{jn}^2 - \sum_{j > ck_n} \zeta_{jn}^2\right| \le  \sum_{j > ck_n} |\eta_{jn}(2\theta_{jn} - \eta_{jn})|\\& \le \left( \sum_{j > ck_n}\eta_{jn}^2\right)^{1/2}\left(2\left( \sum_{j > ck_n} \theta_{jn}^2\right)^{1/2} + \left( \sum_{j > ck_n}\zeta_{jn}^2\right)^{1/2}\right)
\le  C\delta^{1/2}n^{-2r}.
\end{split}
\end{equation}
By (\ref{uh11}) and (\ref{dub2}), we get
\begin{equation}\label{dub3}
 \sum_{j < ck_n} \theta_{jn}^2 
 \ge  \sum_{j < ck_n} \eta_{jn}^2 -  \sum_{j > ck_n} \eta_{jn}^2 -J_n 
 \ge  \sum_{j < ck_n} \eta_{jn}^2 - \delta n^{-2r} - C\delta^{1/2}n^{-2r}.
\end{equation}
Hence, by (\ref{u1}) and Lemma \ref{ld1}, we have
\begin{equation}\label{dop1}\begin{split}&
A_n(\theta) = n^2 \sum_{j=1}^{\infty} \kappa^2_{jn}\theta_{jn}^2 \ge c_3 n k_n^{-1/2}\sum_{j=1}^{ck_n} \theta_{jn}^2\\& \ge  c_3 n k_n^{-1/2}\left(\sum_{j=1}^{ck_n} \eta_{jn}^2 - C\delta^{1/2}n^{-2r}\right) \asymp n k_n^{-1/2} n^{-2r} \asymp 1.
\end{split}
\end{equation}
By Theorem \ref{tq2}, (\ref{dop1}) implies Lemma \ref{ld4}.
\begin{lemma}\label{ld5} Let $\| f_n \| < C n^{-r}$ and let (\ref{ma2}) hold. Then sequence $f_n$ is $n^{-r}$ -- inconsistent.
\end{lemma}
Proof. Let $f_n = \sum_{j = 1}^\infty \theta_{jn} \phi_j$.
Denote
$
 f_{1n} = \sum_{j = 1}^{ck_n} \theta_{jn} \phi_j$.

By Lemma \ref{ld3}, $f_{1n} \in \gamma U$ for some $\gamma  > 0$. If $\|  f_{1n} \| > cn^{-r}$, then, by {\sl i.} in definition of maxiset, $ f_{1n}$ is consistent. Therefore, by Theorem \ref{tq2}, sequence $f_n$ is consistent as well.

Suppose
$\|  f_{1n}  \| = o(n^{-r})$. Then we have
\begin{equation}\label{d7}
n^{2}  \sum_{j=1}^{\infty} \kappa^2_{jn}\theta_{jn}^2 = n^{2}  \sum_{j>ck_n} \kappa^2_{jn}\theta_{jn}^2 + o(1).
\end{equation}
By A1, we have
\begin{equation} \label{d8}
n^{2}  \sum_{j>ck_n} \kappa^2_{jn}\theta_{jn}^2 \le n^2 \kappa^2_{[ck_n],n} \sum_{j>ck_n} \theta_{jn}^2  = o(1)
\end{equation}
as $c \to \infty$ and $n \to \infty$.

By Theorem \ref{tq2}, (\ref{d7}) and (\ref{d8}) imply Lemma \ref{ld5}.
\begin{lemma}\label{ld6} Let sequence $f_{n}$, $cn^{-r}\le \|f_{n}\| \le Cn^{-r}
$, be consistent. Then (\ref{ma1}) holds.
\end{lemma}
Proof. Suppose that, for subsequence $f_{n_i}$, (\ref{ma1}) does not valid. Define sequence $k_{n_i} \asymp n_i^{2-4r}$. Let $f_{1n_i} = \sum_{j=1}^{k_{n_i}} \theta_{jn_i}\,\phi_j$.

If
\begin{equation} \label{d9}
 \sum_{j=1}^{k_{n_i}} \theta_{jn_i}^2 \asymp k_{n_i}^{-2s} \asymp n_i^{-2r},
 \end{equation}
 then, by Lemma \ref{ld3} and {\sl i.} in definition of maxiset, the sequence $f_{n_i}$ is consistent and (\ref{ma1}) holds.

 If (\ref{d9}) does not hold, then, implementing estimates (\ref{d7}), (\ref{d8}) and Theorem \ref{tq2}, we get that sequence $f_n$ is inconsistent.
\begin{lemma}\label{ld7} Let sequence $f_{n}$, $cn^{-r}\le \|f_{n}\| \le Cn^{-r},
$ be inconsistent. Then (\ref{ma2}) holds and sequence $f_n + f_{1n}$ is consistent for any sequence $f_{1n} \in \gamma U, \gamma > 0$.
\end{lemma}
Proof.  Let $f_{1n} = \sum_{j = 1}^{\infty} \eta_{jn} \phi_j \in \gamma U$. If $f_n$ is inconsistent,  then, by Lemma \ref{ld2}, for $k_{n} \asymp n^{2-4r}$ and any $c$, we have $\|\bar f_n \| = o(n^{-r})$ with $\bar f_n = \sum_{j = 1}^{ck_n} \theta_{jn} \phi_j$. By Lemma \ref{ld1} $\sum_{j=1}^{c_2 k_n} \eta_{jn}^2  > c_3^2 n^{-2r}$ for some $c_2, c_3$.

Hence, we get
\begin{equation} \label{d109}
\|\bar f_n + \bar f_{1n}\| = \|\bar f_{1n}\|(1 + o(1)).
\end{equation}
By (\ref{uh11}), for any $\delta > 0$,  for any $cU$, there is $c_1$ such that there holds $\| \tilde f_{1n}  \| < \delta n^{-r}$ where $\tilde f_{1n} = \sum_{j=k_n}^\infty \eta_{jn} \phi_j$ with $k_n = [c_1 n^{2-4r}]$.  Hence
\begin{equation} \label{d110}
|\,\|\tilde f_n + \tilde f_{1n}\| - \|\tilde f_{1n}\|\,| \le \|\tilde f_n\| \le \delta n^{-r}.
\end{equation}
Now (\ref{d109}) and (\ref{d110}) implies (\ref{ma2})

By (\ref{d109}), we have
\begin{equation} \label{d111}
\sum_{j=1}^\infty \kappa_{jn}^2 (\theta_{jn} + \eta_{jn})^2 \ge \sum_{j=1}^{k_n} \kappa_{jn}^2 (\theta_{jn} + \eta_{jn})^2 \asymp \kappa^2 \|\bar f_n + \bar f_{1n}\|^2 \asymp \kappa^2 \| \bar f_{1n}\|^2 \asymp \kappa^2 n^{-2r}.
\end{equation}
 By Theorem  \ref{tq2}, this implies that sequence of alternatives $f_n + f_{1n}$ is consistent.

{\sl Proof of Theorem \ref{tq4}}. For any $\delta >0, c_1, c$ and $C$, we can choose $\gamma$ such that, for any sequence of alternatives $f_n =\sum_{j=1}^\infty \theta_{jn} \phi_j$, $cn^{-r}\le \|f_{n}\| \le Cn^{-r}$, we have  $f_{n\gamma} =\sum_{j=1}^\infty \eta_{jn} \phi_j$  with $\eta_{jn} = \theta_{jn}$ for $1 \le j < ck_n$ and
\begin{equation} \label{qr2}
\sum_{j>c_1k_n} \eta_{jn}^2 \le \delta.
\end{equation}
This allows   to prove {\sl iii.} using the same reasoning as in the proof of Theorem \ref{tq3}.

By (\ref{d8}), for any $\delta > 0$ there is $c$ such that
\begin{equation} \label{tr2}
n^{2}  \sum_{j>ck_n} \kappa^2_{jn}\theta_{jn}^2 \le \delta.
\end{equation}
By (\ref{qr2}), we get also
\begin{equation} \label{sr2}
n^{2}  \sum_{j>ck_n} \kappa^2_{jn}\eta_{jn}^2 \le \delta.
\end{equation}
Hence, by  Theorem \ref{tq2}, implementing  estimates similar to the proof of Lemma \ref{ld3}, we get Theorem \ref{tq4}.
  \subsection{Proof of Theorems of section 5}
  Denote
$$
T_{1n}(f) = \int_0^1\Bigl(\frac{1}{h_n}\int K\Bigl(\frac{t-s}{h_n}\Bigr)f(s) ds\Bigr)^2 dt.
$$
Define the set
$$
Q_{nh_n} = \{f: T_{1n}(f) > \rho_n, f \in  L_2^{per}(R^1)\}.
$$
Proof of Theorem \ref{tk1} is based on the following Theorem \ref{tk2} on asymptotic minimaxity of kernel-based tests \cite{er03}.
\begin{thm}\label{tk2} Let $h_n^{-1/2}n^{-1} \to 0$, $h_n \to 0$ as $n \to \infty$.
Let
\begin{equation}\label{}
0 < \liminf_{n \to \infty}  n \rho_n h_n^{1/2} \le \limsup_{n \to \infty} n\rho_nh_n^{1/2} < \infty.
\end{equation}
Then the family  of kernel-based tests $L_n = \chi\{T_n(Y_n) \ge x_\alpha\}, \alpha(L_n) = \alpha(1 + o(1))$, is asymptotically minimax for the sets of alternatives $Q_{nh_n}$.

There holds
\begin{equation}\label{}
\beta(L_n,Q_{nh_n})  = \Phi(x_\alpha - \kappa^{-1} \sigma^{-2}n h_n^{1/2}\rho_n)(1 + o(1)).
\end{equation}
Here $x_\alpha$ is defined the equation $\alpha = 1 - \Phi(x_\alpha)$.

Moreover, for each $f_n \in L_2^{per}(R^1)$,  there holds
\begin{equation}\label{33}
\beta(L_n,f_n)  = \Phi(x_\alpha - \kappa^{-1} \sigma^{-2}n h_n^{1/2}\rho_n)(1 + o(1))
\end{equation}
uniformly on $f_n$ such that $T_{1n}(f_n) = \rho_n(1 + o(1))$.
\end{thm}

{\sl Proof of Theorem \ref{tk1}. Sufficiency.} Let $f_n \in \mathbb{B}^s_{2\infty}$ and let $\|f_n\| \asymp n^{-r}$.
By Theorem  \ref{tk2}, the consistency takes place if
\begin{equation}\label{kp2}
\rho_n \asymp  n^{-1}h_n^{-1/2} \asymp n^{-2r}.
\end{equation}

 We shall explore the problem in terms of sequence model.

For $-\infty < j < \infty$, denote
$$
\hat K(jh) =  \frac{1}{h}\int_{-1}^1 \exp\{2\pi ijt\} K\left(\frac{t}{h}\right) dt,
$$
$$
y_j = \int_0^1 \exp\{2\pi ijt\} dY_n(t),
$$
$$
\xi_j = \int_0^1 \exp\{2\pi ijt\} dw(t),
$$
$$
\theta_j = \int_0^1 \exp\{2\pi ijt\} f(t) dt.
$$
Denote $Y_n = \{y_j\}_{-\infty}^\infty$.

In this notation we can write our sequence model in the following form
\begin{equation}\label{}
y_j = \hat K(jh) \theta_j + \sigma n^{-1/2} \hat K(jh) \xi_j, \quad -\infty \le j < \infty.
\end{equation}
and
$$
T_n(Y_n) = nh_n^{1/2}\sigma^{-2}\kappa^{-1} \Bigl(\sum_{j=-\infty}^\infty |\hat K^2(jh) y_j^2|  -  n^{-1}\sigma^2 \sum_{j=-\infty}^\infty |\hat K^2(jh)|\Bigr).
$$
The function $\hat K(\omega)$, $\omega \in R^1,$ is analytic  and $\hat K(0) = 1$. Therefore there is an interval $(-b,b), 0 <b< \infty,$ such that $\hat K(\omega) \ne 0$ for all $\omega  \in (-b,b)$.

We have
\begin{equation}\label{}
\sum_{|j| > bh_n^{-1}} |\theta_j|^2 =O(b^{-2s}h_n^{2s})
\end{equation}
Therefore, there exists $c>0$ such that, for $h_n < bcn^{-2/(1+4s)}$, there holds
\begin{equation}\label{ke1}
\rho_n \asymp n^{-2r} \asymp \sum_{|j| < b h_n^{-1}}\,|\theta_j|^2 \asymp \sum_{|j| < bh_n^{-1}}\,|\hat K(jh_n)\,\theta_j|^2 \asymp n^{-1}h_n^{1/2}.
\end{equation}
By (\ref{33}) and (\ref{ke1}), we get sufficiency.

{\sl Proof of necessary conditions}.
Suppose the opposite. Then there are vector $\theta = \{\theta_j\}_{j=1}^\infty$  and a sequence $m_l$, $m_l \to \infty$ as $l \to \infty$, such that
\begin{equation}\label{}
m_l^{2s} \sum_{|j|\ge m_l}^\infty |\theta_j|^2 = C_l
\end{equation}
with $C_l \to \infty$ as $l \to \infty$.

It is clear that  we can define a sequence $m_l$ such that
\begin{equation}\label{}
m_l^{2s} \sum_{m_l \le |j| \le 2m_{l}} |\theta_j|^2 > \delta C_l
\end{equation}
where $\delta >0$  does not depend on $l$ .

Define a sequence $\eta_l = \{\eta_{jl}\}_{j=-\infty}^\infty$ such that
 $\eta_{jl} = \theta_j,  |j| \ge m_{l}$, and $\eta_{jl} = 0$ otherwise.

Denote
$$
\tilde f_l(x) = f_{l}(x,\eta_l) = \sum_{j=-\infty}^\infty \eta_{jl} \exp\{2\pi ijx\}.
$$
 For alternatives $\eta_l$ we define $n_l$ such that $\|\eta_l\| \asymp n_l^{-r}$.

 Then
\begin{equation}\label{k100}
n_l \asymp C_l^{-1/(2r)} m_l^{s/r}
\end{equation}
We have $|\hat K(\omega)| \le \hat K(0) = 1$ for all $\omega \in R^1$ and $|\hat K(\omega)| > c > 0$ for $ |\omega| < b$. Hence, if we put $h_l= h_{n_l} =2^{-1}b^{-1}m_l^{-1}$, then there is $C > 0$  such that, for all $h> 0$, there holds
\begin{equation} \label{}
T_{1n_l}(\tilde f_l,h_l) = \sum_{j=-\infty}^\infty |\hat K(jh_l)\,\eta_{jl}|^2 > C \sum_{j=-\infty}^\infty |\hat K(jh)\,\eta_{jl}|^2 = C T_{1n_l}(\tilde f_l,h).
\end{equation}
Thus we can choose $h = h_l$ for further reasoning.

We have
\begin{equation}\label{k101}
\rho_{n_l} =  \sum_{|j|>m_l}\, |\hat K(jh_l) \eta_{jl}|^2  \asymp \sum_{j=m_l}^{2m_l} |\eta_{jl}|^2 \asymp n_l^{-2r}.
\end{equation}
If we put in estimates (\ref{u8})--(\ref{u10}), $k_l = [h_{n_l}^{-1}]$ and $k_l = m_l$, then we get
\begin{equation}\label{k102}
h_{n_l}^{1/2} \asymp C_l^{(2r-1)/2}n^{2r-1}.
\end{equation}
By (\ref{k101}) and (\ref{k102}), we get
\begin{equation}\label{}
n_l\rho_{n_l}h_{n_l}^{1/2} \asymp C_l^{-(1-2r)/2}.
\end{equation}
By Theorem \ref{tk2}, this implies inconsistency of hypothesis and alternatives $\eta_l$.

{\sl Proof of Theorem \ref{tk3}}. Test statistics $T_n(Y_n)$ are quadratic forms. Therefore, for the proof of {\sl i.} and {\sl ii.}, we can implement the same reasoning as in the proof of Theorem \ref{tq1}. Theorem \ref{tk2} can be treated as a version of Theorem \ref{tq2} with $\kappa_{jn}^2 = |\hat K(jh_n)|^2$ and $k_n = [h_n^{-1}]$.

Since it is known only that $|\hat K(\omega)| > c > 0$ for $|\omega| < b$, we are forced to make small differences in the reasoning.
The differences are the following. In version of Lemma \ref{ld2} and in the proof of version of Lemma \ref{ld4} we need to suppose additionally that $c< b$. In the proof of Lemma \ref{ld5} one needs to replace $\kappa^2_{ck_n,n}$ with $\sup_{|\omega| > c}| \hat K(\omega)|^2 h_n^{1/2}$.

Proof of Theorem \ref{tk4} is akin to  the proof of Theorem \ref{tq4} and is omitted.
    \subsection{Proof of Theorems of section 6}
  Sufficiency  in  Theorem \ref{chi1} has  been proved  Ingster \cite{ing87}.

The proof of necessary condition in  Theorem \ref{chi1} will be based on Theorem \ref{chi2} provided below. Theorem \ref{chi2} is a summary of results of Theorems 2.1 and 2.4 in Ermakov \cite{er97}.

Denote  $p_{in} =  F(i/k_n) -  F((i-1)/k_n), 1 \le i \le k_n$.

Define the sets of alternatives
$$
Q_n(b_n) = \Bigl\{ F: T_n(F)= n k_n  \sum_{i=1}^{k_n} (p_{in}  - 1/k_n)^2 \ge b_n\Bigr\}.
$$
The definition of asymptotic minimaxity of test is the same as in section 3.

Define the tests
$$
K_n = \chi(2^{-1/2}k_n^{-1/2}(T_n(\hat F_n) - k_n + 1) >x_\alpha)
$$
where $x_\alpha$ is defined the equation $\alpha = 1- \Phi(x_\alpha)$.
\begin{thm}\label{chi2} Let $k_n^{-1} n^2 \to \infty$ as $n \to \infty$. Let
\begin{equation}\label{}
0 < \liminf_{n \to \infty} k_n^{-1/2}b_n \le \limsup_{n \to \infty} k_n^{-1/2} b_n < \infty.
\end{equation}
Then $\chi^2$-tests  $K_n$ are asymptotically minimax for the sets of alternatives $Q_n(b_n)$.

There holds
\begin{equation}\label{}
\beta(K_n,F) = \Phi(x_\alpha - 2^{-1/2}k_n^{-1/2}T_n(F))(1 +o(1))
\end{equation}
uniformly in $F$ such that $ck_n^{1/2} \le T_n(F) \le C k_n^{1/2}$.
\end{thm}
 For any complex number $a = b+ id$ denote $\bar a = b - id$.

We have
\begin{equation}\label{}
n^{-1}k^{-1}T_n(F) = \sum_{l=0}^{k-1}\left(\int_{l/k}^{(l+1)/k} f(x) dx \right)^2.
\end{equation}
We can write $f(x)$ in terms of Fourier coefficients
\begin{equation}\label{}
f(x) = \sum_{j=-\infty}^\infty \theta_j \exp\{2\pi ijx\}.
\end{equation}
Then
\begin{equation}\label{}
\int_{l/k}^{(l+1)/k} f(x) dx = \sum_{j=-\infty}^\infty \frac{\theta_j}{2\pi ij}\exp\{2\pi ijl/k\} ( \exp\{2\pi ij/k\} - 1).
\end{equation}
Hence
\begin{equation}\label{}
\begin{split}&
n^{-1}k_n^{-1}T_n(F) =\sum_{l=0}^{k_n-1}\Bigl(\sum_{j \ne 0}\frac{\theta_j}{2\pi ij}\exp\{2\pi ijl/k_n\} ( \exp\{2\pi ij/k_n\} - 1\}\Bigr)\\& \times
\Bigl(\sum_{j \ne 0} \frac{-\bar\theta_j}{2\pi ij}\exp\{-2\pi ijl/k_n\} \exp\{-2\pi ij/k_n\} - 1\Bigr) = J_1 + J_2,
\end{split}
\end{equation}
with
\begin{equation}\label{}
\begin{split}&
J_1=\sum_{l=0}^{k_n-1}\sum_{m=-\infty}^\infty\,\,\, \sum_{{j \ne mk_n}, {j_1 = j -mk_n}} \frac{\theta_j \bar\theta_{j_1}}{4\pi^2jj_1}\exp\{2\pi ilm\}\\&\times(\exp\{2\pi ij/k_n\} - 1)(\exp\{-2\pi ij_1/k_n\} - 1)\\&= k_n \sum_{m=-\infty}^\infty \sum_{j \ne mk_n} \frac{\theta_j \bar\theta_{j-mk_n}}{4\pi^2j(j-mk_n)}(2 - 2 \cos(2\pi j/k_n))
\end{split}
\end{equation}
and
\begin{equation}\label{ux}
\begin{split}&
J_2=\sum_{l=0}^{k_n-1}\sum_{j=-\infty}^\infty \sum_{j_1 \ne j-m k_n} \frac{\theta_j \bar\theta_{j_1}}{4\pi^2 j j_1} \exp\{2\pi i(j - j_1)l/k_n\}\\&\times(\exp\{2\pi ij/k_n\} - 1)(\exp\{-2\pi ij_1/k_n\} - 1) =0,
\end{split}
\end{equation}
where $\sum_{j_1 \ne j-m k_n}$ denotes summation over all $j_1$ such that $j -j_1 \ne m k_n$ for all integers $m$.
In the last equality of  (\ref{ux}) we make use of the identity
\begin{equation}\label{}
 \sum_{l=0}^{k-1}\exp\{2\pi i(j - j_1)l/k\} = \frac{\exp\{2\pi i(j - j_1) k/k\} -1}{exp\{2\pi i(j - j_1)/k\} -1} =0,
\end{equation}
if $j-j_1 \ne mk, -\infty < m < \infty$.

For any c.d.f $F$ denote $\tilde F_k$ c.d.f. with the density
$$
 1 + \tilde f_k(x)  = 1 + \sum_{|j| > k} \theta_j \exp\{2\pi ijx\}.
$$
Suppose the opposite. Then there is sequence $i_l, i_l \to \infty$ as $l \to \infty$, such that
\begin{equation}\label{}
i_l^{2s} \|\tilde f_{i_l}\|^2  = C_l,
\end{equation}
with $C_l \to \infty$ as $l \to \infty$.

By  Theorem \ref{chi2}, it suffices to show that $k_{l}^{-1/2}T_{n_l}(\tilde F_{i_l}) = o(1)$ with $n_l$ defined the equation
\begin{equation}\label{}
\|\tilde f_{i_l}\|^2= \sum_{|j| > i_l}|\theta_j|^2\asymp n_l^{-2r}
\end{equation}
and $k_l = k_{n_l}\asymp n_l^{2- 4r}$.

Then $i_l \asymp C_l^{\frac{1}{2s}} n_l^{2-4r}$. Denote $d_l = i_l/k_l$.

Denote $\eta_j = \theta_j$ if $|j|> i_l$ and $\eta_j =0$ if $|j| < i_l$.

We have $|j - mk_l| \ge |m-1|  k_l$ if $j \ge 2k_l$, $m \ne 1$ or $j < 0$, $m \ne 1$ or $j \le -2k_l$, $m \ne -1$ or $j > 0, m \ne -1$. We also  have $\eta_{j-mk_l}/(j -mk_l) = 0$ if $0 < j < 2k_l$, $m=1$  or $-2k_l < j < 0$, $m = -1$.

Hence, implementing  $|j - mk_l| \ge |m-1|  k_l$ in the first and in the third inequalities of (\ref{hu}), we get

\begin{equation}\label{hu}
\begin{split}&
n_l^{-1}k_{n_l}^{-1}T_{n_l}(\tilde F_{i_l}) =\sum_{m=-\infty}^\infty \,\,\sum_{j \ne mk_l, |j|>k_l} \frac{\eta_j \bar\eta_{j-mk_l}}{4\pi^2j(j-mk_l)}(2 -2 \cos(2\pi j/k_l))\\&
\le C k_l^{-1} \sum_{j=-\infty}^\infty \Bigl|\frac{\eta_j}{j}\Bigr| \sum_{m \ne 0} \Bigl|\frac{\eta_{j+mk_l}}{m}\Bigr|
\\&
\le C k_l^{-1} \sum_{j=1}^{k_l} \sum_{m_1=-\infty}^\infty \Bigl|\frac{\eta_{j+m_1k_l}}{j+m_1k_l}\Bigr|  \sum_{m + m_1 \ne 0}
\Bigl|\frac{\eta_{j+(m+m_1)k_l}}{m+m_1}\Bigr| \\&
\le Ck_l^{-2} \sum_{j=1}^{k_l} \sum_{|m_1| > d_l}\Bigl|\frac{\eta_{j+m_1k_l}}{m_1}\Bigr|  \sum_{|m + m_1| >  d_l}
\Bigl|\frac{\eta_{j+(m+m_1)k_l}}{m+m_1}\Bigr|
\\& \le Ck_l^{-2} \sum_{j=1}^{k_l}  \Bigl( \sum_{|m_1| > d_l} |\eta_{j+m_1k_l}|^2  \Bigl(\sum_{|m + m_1| >  d_l} \Bigl|\frac{\eta_{j+(m+m_1)k_l}}{m+m_1}\Bigr|\Bigr)^2\Bigr)^{1/2} \Bigl(\sum_{|m_1| > d_l} m_1^{-2}\Bigr)^{1/2}
\\&
 \le Ck_l^{-2} \sum_{j=1}^{k_l} \Bigl(\sum_{|m_1| > d_l}|\eta_{j+m_1k_l}|^2  \sum_{|m + m_1| >  d_l} |\eta_{j+(m+m_1)k_l}|^2 \sum_{|m_2 + m_1| >  d_l}(m_2+m_1)^{-2}\Bigr)^{1/2}\\&\times
 \Bigl( \sum_{|m_1| > d_l} m_1^{-2}\Bigr)^{1/2}\\& \le Ck_l^{-2} \sum_{j=1}^{k_l} \sum_{|m_1| > d_l}|\eta_{j+m_1k_l}|^2 \sum_{|m| > d_l} m^{-2}
 \le C k_l^{-1} i_l^{-1}\sum_{|j| > i_l} |\theta_j|^2.
\end{split}
\end{equation}
Hence
\begin{equation}\label{}
k_{l}^{-1/2}T_{n_l}(\tilde F_{i_l}) \le k_{l}^{1/2} i_l^{-1} n_l \sum_{|j| > i_l} |\theta_j|^2 \asymp k_{l}^{1/2} i_l^{-1} n_l^{1-2r} \asymp C_l^{-1/2s}.
\end{equation}
By  Theorem \ref{chi2}, this implies the necessary conditions.

{\sl Proof of Theorem \ref{chi3}}. Analysis of the proof of Lemmas \ref{ld1} - \ref{ld7} shows that, for the proof of Theorem \ref{chi3} it suffices  to prove Lemmas \ref{lchi3} - \ref{lchi2} provided below.

Let
$$ f_n = \sum_{j=-\infty}^\infty \theta_{jn} \phi_j, \quad \phi_j(x) = \exp\{2\pi i\,j\,x\,\}, \quad x \in (0,1).
$$
Let $k_n = \Bigl[n^\frac{2}{1 + 4s}\Bigl]$.
\begin{lemma}\label{lchi3} If sequence $f_n$, $cn^{-r} < \|f_n\| < C n^{-r}$, is consistent, then there are $c_1$ and  $c_2$ such that there holds
\begin{equation}\label{h1}
\sum_{|j| < c_2 k_n} |\theta_{jn}|^2 > c_1 n^{-r}.
\end{equation}
\end{lemma}
\begin{lemma}\label{lchi1} If, for the sequence $f_n$,  there are $c_1$ and  $c_2$ such that (\ref{h1}) holds, then there is sequence $k_n \asymp n^{2-4r}$ such that $f_n$ is consistent for the chi-squared test statistics $T_n$.
\end{lemma}
If (\ref{h1}) holds, one can put $f_{1n} = \sum_{|j| < c_2 k_n} \theta_{jn} \phi_j$ in (\ref{ma1}) and {\sl i.} will be hold.

If (\ref{ma1}) holds, then (\ref{uh11}) implies that (\ref{h1}) holds. Therefore, by Lemma \ref{lchi3} sequence of alternatives $f_n$ is consistent
\begin{lemma}\label{lchi2} Sequence $f_n$, $cn^{-r} < \|f_n\| < C n^{-r}$, is inconsistent, iff, for all $c_2$, there holds
\begin{equation}\label{h2}
\sum_{|j| < c_2 k_n} |\theta_{jn}|^2 = o (n^{-r})
\end{equation}
as $n \to \infty$.
\end{lemma}
If (\ref{h2}) holds, then, arguing similarly to the proof of Lemma \ref{ld7}, we get that (\ref{ma2}) holds and sequence of alternatives is consistent.

Suppose (\ref{ma2}) holds and $\| f_n\| < C n^{-r}$. Denote $f_{1n} = \sum_{j=1}^{c_0 k_n} \theta_{jn} \phi_j$. If $\| f_{1n}\| > c_1  n^{-r}$, then sequence $f_n$ is consistent by (\ref{h1}).

If $\| f_{1n}\| = o(n^{-r})$ for each $c_0$, then sequence $f_n$ is is inconsistent by (\ref{h2}).
This implies necessary conditions in {\sl ii.}

{\sl Proof of Lemma \ref{lchi1}}. Let (\ref{h1}) hold.
For any $a >  0$ denote
$$  \tilde f_{n,ak_n}  = \sum_{|j| > a k_n} \theta_{jn}\phi_j$$
and denote
$$
f_{n,c_1k_n,C_1k_n}=\tilde f_{n,c_1k_n} - \tilde f_{n,C_1k_n}, \quad \bar f_n =\bar f_{n,c_1k_n} = f_n - \tilde f_{n,c_1k_n},
$$
with $C_1 > c_1$.

Let $T_n$ be the chi-squared test statistics with the number of cells $l_n= [c_3 k_n]$  with $c_1 < c_3 < C_1$.

We have
\begin{equation}\label{h3}
T_n^{1/2}(\bar f_{n}) - T_n^{1/2}(f_{n,C_1k_n,c_1k_n})- T_n^{1/2}(\tilde f_{n,C_1k_n}) \le T_n^{1/2}(f_{n}).
\end{equation}
Denote
$$
\bar p_{jn} = \frac{1}{k_n} \int_{(j-1)/k_n}^{j/k_n} \bar f_n(x) dx.
$$
By Lemmas 3 and 4 in section 7 of Ulyanov \cite{ul}, we have
\begin{equation}\label{h5}
S_n(\bar f_n)  \doteq \sum_{j=1}^{k_n} \int_{(j-1)/k_n}^{j/k_n}(\bar f_n(x) - \bar p_{jn})^2 \, dx \le 2 \omega^2\Bigl(\frac{1}{k_n}, \bar f_n\Bigr).
\end{equation}
Here
$$\omega^2(h,f) = \int (f(t+h) - f(t))^2\,dt, \quad h>0,$$
for any $f \in L_2^{per}$.
If $f = \sum_{j=-\infty}^\infty \theta_j \phi_j,
$ then
\begin{equation}\label{h6}
\omega^2(h,f) =  2\sum_{j=1}^\infty |\theta_j|^2\,(2  - 2\cos (jh)).
\end{equation}
 Since $1 - cos( x) \le x^2$, then, by (\ref{h5}) and (\ref{h6}), we have
\begin{equation}\label{h8}
\|\bar f_n\| - k_n^{-1/2} n^{-1/2}T_n^{1/2}(\bar f_n) \le S_n^{1/2}(\bar f_n) \le c_1 c_3^{-1} \|\bar f_n\|.
\end{equation}
By (\ref{hu}), we get
\begin{equation}\label{h10}
k_n^{-1} n^{-1}T_n(\tilde f_n) < C_1^{-1} c_3\|\tilde f_n\|^2 < C C_1^{-1} c_3 n^{-2r}.
\end{equation}
 We have
 \begin{equation}\label{h7}
k_n^{-1/2}n^{-1/2}T_n^{1/2}( f_{n,c_1k_n,C_1k_n}) \le  \|f_{n,c_1k_n,C_1k_n}\|.
\end{equation}
 Fix $\delta, 0< \delta <1$ and fix $c_2$. There are at most $2[\delta^{-1}]$ intervals $[c_2\delta^{-2i},C_2\delta^{-2i-2}]$, $0 \le i \le 2\delta^{-1}$ such that for one of them, for $c_1 = c_2\delta^{-2i}$ and $C_1 = C_2\delta^{-2i-2}$ there holds
\begin{equation}\label{h9}
\sum_{c_1k_n < |j| < C_1 k_n} |\theta_{jn}|^2 =  \|f_{n,c_1k_n,C_1k_n}\|^2 < C\delta n^{-r}.
\end{equation}
For any $c_1$ and $C_1$ such that (\ref{h9}) holds, we put $c_3 = c_1 \delta$.

Since the choice of $\delta$ was arbitrary, then, by (\ref{h3}), (\ref{h8}),(\ref{h10}) and (\ref{h9}) together, we get $k_n^{-1/2}T_n(f_n) \asymp 1$. By Theorem \ref{chi2}, this implies sufficiency.

{\sl Proof of Lemma \ref{lchi2}. Sufficiency.} In the proof of sufficiency we choose test statistics $T_n$ with sufficiently large number of cells $k_n$. It is clear that we can always make additional partitions of cells and test statistics with these  additional partitions of cells will be also consistent if the number of cells will have the same order $n^{2-4r}$.

 We have
\begin{equation}\label{h11}
 T_n^{1/2}(f_{n}) \le T_n^{1/2}(\bar f_{n,Ck_n}) + T_n^{1/2}(\tilde f_{n,Ck_n}).
\end{equation}
By (\ref{hu}), we have
\begin{equation}\label{h12}
n^{-1}k_n^{-1} T_n(\tilde f_{n,Ck_n}) \le C^{-1} \|\tilde f_n \|^2  \le C^{-1} n^{-2r}.
\end{equation}
We have
\begin{equation}\label{h13}\|\bar f_{n,Ck_n}\|^2 \ge n^{-1}k_n^{-1} T_n^{1/2}(\bar f_{n,Ck_n}).
\end{equation}
By Theorem \ref{chi2}, (\ref{h2}) and  (\ref{h11}) - (\ref{h13}) together implies
inconsistency of sequence $f_n$.

If $f_n$ is inconsistent, then (\ref{h2}) follows from  Lemma \ref{lchi1}.

 Lemma \ref{lchi3} follows from sufficiency statement of Lemma \ref{lchi2}.
   \subsection{Proof of Theorems of section 7}
   {\sl Proof of Theorem \ref{tom1}}.
We can write the functional $T^2(F-F_0)$ in the following form (see Ch.5,  Shorack and Wellner \cite{wel})
\begin{equation}\label{om1}
T^2(F-F_0) = \int_0^1\int_0^1(\min\{s,t\} - st) f(t) f(s) ds dt
\end{equation}
with $f(t) = d(F(t) - F_0(t))/dt$.

If we consider the expansion of function
\begin{equation}\label{om22}
f(t) = \sqrt{2} \sum_{j=1}^\infty \theta_j \cos(\pi j t), \quad \theta = \{\theta_j\}_{j=1}^\infty
\end{equation}
on eigenfunctions of operator with the kernel $\min\{s,t\} -st$, then we get
\begin{equation}\label{om3}
nT^2(F-F_0) = n\sum_{j=1}^\infty\frac{\theta_j^2}{\pi^2 j^2}
\end{equation}
Proof of {\sl i.} For this setup {\sl i.} has the following form
\vskip 0.3cm
{\sl i.} for all $\theta \in U, ||\theta|| > n^{-r}$, there holds
\begin{equation}\label{omx}
n\sum_{j=1}^\infty\frac{\theta_j^2}{\pi^2 j^2} > c,
\end{equation}
Note that (\ref{omx}) can be replaced with the following condition
\begin{equation}\label{om5}
n \sum_{k=1}^\infty 2^{-2k}\sum_{j = 2^k+1}^{2^{k+1}} \theta_j^2 >c
\end{equation}
and we suppose that
\begin{equation}\label{om6}
\sum_{k=1}^\infty \sum_{j = 2^k+1}^{2^{k+1}} \theta_j^2 > n^{-2r}.
\end{equation}
and
\begin{equation}\label{om7}
2^{2ls}\sum_{k=l}^\infty \sum_{j = 2^k+1}^{2^{k+1}} \theta_j^2 \le P_0
\end{equation}
for all $l$.

Denoting $\beta_k^2=\sum_{j = 2^k+1}^{2^{k+1}} \theta_j^2$ we can rewrite (\ref{om5})-(\ref{om7}) in the following form
\begin{equation}\label{om8}
n\sum_{k=1}^\infty 2^{-2k}\beta_k^2 > c
\end{equation}
and we suppose that
\begin{equation}\label{om9}
\sum_{k=1}^\infty \beta_k^2 > n^{-2r}
\end{equation}
and
\begin{equation}\label{om10}
f = \{\beta_j\}_{j=1}^\infty \in W = \Bigl\{f : \sup_l 2^{2ls}\sum_{j=l}^\infty \beta_j^2  \le P_0, f = \{\beta_j\}_{j=1}^\infty \Bigr\}.
\end{equation}
The infimum of the left-hand side of (\ref{om8}) is attained for $\beta =\{\beta_k\}_{k=1}^\infty$ such that, for some $k=k_0$ there hold $P_0/2 < 2^{2k_0s} \beta_{k_0}^2 \le P_0$ and $\beta_k =0$ for $k <k_0$.

Hence,  by (\ref{om8}), we get
\begin{equation}\label{om11}
\beta_{k_0}^2 \asymp 2^{-2k_0s}P_0\asymp n^{-2r}.
\end{equation}
Therefore
\begin{equation}\label{om12}
2^{2k_0} \asymp  n^{2r/s} \asymp n^{1-2r}.
\end{equation}
Hence we get
\begin{equation}\label{om13}
n\sum_{k=1}^\infty 2^{-2k}\beta_k^2  \asymp n2^{-2k_0}\beta_{k_0}^2 \asymp n 2^{-2k_0} n^{-2r} \asymp 1.
\end{equation}
This implies {\sl i.}
\vskip 0.3cm
{\sl Proof of necessary conditions}.
Suppose the opposite. Then there is a sequence $m_i$ such that
\begin{equation}\label{om14}
2^{2m_i s} \sum_{k=m_i}^\infty \beta_k^2  = C_i \to  \infty
\end{equation}
as $i  \to \infty$.

Define sequence $n_i$ such that
\begin{equation}\label{om17}
n_i^{-2r} \asymp  \sum_{k=m_i}^\infty \beta_k^2 \asymp C_{m_i} 2^{-2m_is}.
\end{equation}
Then
\begin{equation}\label{om18}
2^{-2m_i} \asymp C_{m_i}^{-1/s}n_i^{-2r/s} \asymp C_{m_i}^{-1/s}n_i^{2r-1}.
\end{equation}
By (\ref{om17}) and (\ref{om18}), we get
\begin{equation}\label{}
n_i\sum_{k=m_i}^\infty 2^{-2k}\beta_k^2 < C n_i 2^{-2m_i} \sum_{k=m_i}^\infty \beta_k^2 \asymp C_{m_i}^{-1/s}.
\end{equation}
This implies necessary condition.

{\sl Proof of Theorem \ref{tom2}}.
It suffices to prove {\sl i.} and {\sl ii.} in terms of $f_n = \{\beta_{jn}\}_{j=1}^\infty$. In this case  {\sl i.}  and {\sl ii.} in definition of perfect maxisets have similar form. The unique difference is that we replace the set $U$ with the set $W$. The proof of {\sl i.}  and {\sl ii.} is based on versions Lemmas \ref{ld1} -- \ref{ld7} adapted for this setup. The statements of these Lemmas  is the same or almost the same as the statement of Lemmas \ref{ld1} -- \ref{ld7}. Their proofs   represents slight modification of proofs of Lemmas \ref{ld1} -- \ref{ld7}.

Denote $m = [\log_2 n]$.

Sequence $f_n = \{\beta_{jn}\}_{j=1}^\infty$, $c2^{-rm} \le \|f_n\| \le C2^{-rm}$, is inconsistent if
\begin{equation}\label{om13a}
2^m \sum_{j=1}^\infty 2^{-2j} \beta^2_{jn} \to 0\quad\mbox{\rm as} \quad n \to \infty.
\end{equation}
\begin{lemma}\label{lom1} Let $c2^{-rm} \le \|f_n\| \le C2^{-rm}$ and let $f_n \in c_1 W$. Then there is $k_n= (1/2 -r)m + O(1)$ such that
\begin{equation}\label{om14a}
\sum_{j=1}^{k_n} \beta_{jn}^2 > c_2  2^{-2rm}.
\end{equation}
\end{lemma}
Proof. We have
\begin{equation}\label{om15}
2^{2sk_n} \sum_{j=k_n}^{\infty} \beta_{jn}^2 = C 2^{2rm} \sum_{j=k_n}^{\infty} \beta_{jn}^2 \le c_1.
\end{equation}
Hence
\begin{equation}\label{om16}
\sum_{j=k_n}^{\infty} \beta_{jn}^2 \le C^{-1}c_1 2^{-2rm}
\end{equation}
and (\ref{om14a}) holds with $c_2 = c/2$ if $C > \frac{c}{2c_1}$.
\begin{lemma}\label{lom2} Let $f_n$ be $n^{-r}$-inconsistent for the test statistics $T_n$ with $k_n = (1/2 - r)m + O(1)$ as $n \to \infty$. Then we have
\begin{equation}\label{om17a}
2^{2sk_n} \sum_{j=1}^{k_n} \beta_{jn}^2 \asymp 2^{2rm} \sum_{j=1}^{k_n} \beta_{jn}^2 = o(1)
\end{equation}
as $n \to \infty$.
\end{lemma}
Proof. We have
\begin{equation}\label{om18a}
o(1) =  2^m \sum_{j=1}^\infty 2^{-2j} \beta^2_{jn} \ge 2^m 2^{-2k_n}\sum_{j=1}^{k_n} 2^{-2j} \beta^2_{jn} \asymp 2^{2rm} \sum_{j=1}^{k_n} \beta_{jn}^2.
\end{equation}
This implies Lemma \ref{lom2}.
\begin{lemma}\label{lom3} Let  $f_n = \{\beta_{jn}\}_{j=1}^\infty$ and let $\beta_{jn} = 0$ for $j > k_n = (1/2 -r)m + O(1)$. Let $\| f_n \| \le C 2^{-rm}$. Then there is $cW$ such that $f_n \in cW$. \end{lemma}
Proof of Lemma \ref{lom3} is akin to the proof of Lemma \ref{ld3} and is omitted.

The following Lemmas \ref{lom4} and \ref{lom5} have  almost the same statements as Lemmas \ref{ld4} and \ref{ld5}.
\begin{lemma}\label{lom4}  Let $f_{1n} \in cW$. Let $c_1n^{-r}\le \|f_{1n}\| \le C_1n^{-r}
$ and let (\ref{ma1}) hold. Then sequence $f_n$ is $n^{-r}$-consistent. \end{lemma}
Proof.  If (\ref{ma1}) hold then orthogonality of $f_{1n}$ and $f_n - f_{1n}$ does not imply orthogonality of corresponding vectors in terms of coordinates $\beta_{jn}$. At the same time, arguing similarly to the  proof of Lemma \ref{ld4} one can show that, for any $\delta > 0$ there is $k_n = (1/2 -r)m + O(1)$ such that
\begin{equation}\label{dub5}
\sum_{j=1}^{k_n} \beta_{jn}^2 \ge \sum_{j=1}^{k_n} \beta_{1jn}^2 - C\delta^{1/2}n^{-2r},
\end{equation}
where $\beta_{1jn}$, $1 \le j < \infty$, are coordinates of $f_{1n}$

Therefore, by Lemma \ref{lom1}, we have
\begin{equation}\label{om19}
2^m \sum_{j=1}^\infty 2^{-2j} \beta^2_{jn} \ge 2^m2^{-2k_n} \Bigl(\sum_{j=1}^{k_n}  \beta^2_{jn}  \delta n^{-2r}\Bigr) \asymp 1.
\end{equation}
This implies Lemma \ref{lom4}.
\begin{lemma}\label{lom5} Let $\| f_n \| < C n^{-r}$ and let (\ref{ma2}) hold. Then sequence $f_n$ is $n^{-r}$ - inconsistent.
\end{lemma}
Proof. Denote $\bar f_n = \{\tau_{jn}\}_{j=1}^\infty$ with $\tau_{jn} = \beta_{jn}$ for $j \le k_n$ and $\tau_{jn} = 0$ for $j > k_n$.

Denote $\tilde f_n = f_n - \bar f_n$.

By Lemma \ref{lom1}, we get that, if $\|\bar f_n \| > c n^{-r}$ then $f_n$ is $n^{-r}$-consistent.

Suppose $\|\bar f_n \| = o(n^{-r})$. Then we have
\begin{equation}\label{om20}
2^m \sum_{j=1}^\infty 2^{-2j} \beta^2_{jn} = 2^m \sum_{j=k_n+C}^\infty 2^{-2j} \beta^2_{jn} + o(1)
\end{equation}
and
\begin{equation}\label{om21}
2^m \sum_{j=k_n+C}^\infty 2^{-2j} \beta^2_{jn} \le 2^{-C} 2^m 2^{-2k_n}\sum_{j=k_n+C}^\infty \beta^2_{jn} = o(1)
\end{equation}
as $C \to \infty$ and $n \to \infty$. This completes proof of Lemma \ref{lom5}.

The statements of versions of Lemmas \ref{ld6} and \ref{ld7} for this setup is the same. Their proofs are also completely follow the same lines. We omit this reasoning.

{\sl Proof of Theorem \ref{tom2b}}. To implement Hungary construction we need some statement on uniform continuity   of limit distributions of statistics $T_n$ if alternatives hold. This statement is provided in the following Lemma \ref{lom10}.

Denote $b(t)$ Brownian bridge, $t \in (0,1)$.
\begin{lemma}\label{lom10} Assume {\rm B1}. Then the densities of $T^2_n(b(F_n(t)) + \sqrt{n}(F_n(t) - t))$ are uniformly bounded.
\end{lemma}
Proof. We have
\begin{equation}\label{pl1}
\begin{split}&
T^2_n(b(F_n(t)) + \sqrt{n}(F_n(t) - t)) = \int_0^1 (b(F_n(t)) + \sqrt{n}(F_n(t) - t))^2\, dt\\&
 = \int_0^1 \left(\sqrt{2}\sum_{k=1}^\infty \xi_k \frac{\sin(\pi k F_n(t))}{k\pi} + n^{1/2}(F_n(t) - t)\right)^2\, dt,
 \end{split}
\end{equation}
where
$\xi_k = \sqrt{2}\int_0^1 b(t) \sin(\pi k t) \, dt$.

Hence, we have
\begin{equation}\label{pl2}
T^2(\xi_1,\xi_2,J_n) = a_n\xi_1^2 + 2b_n\xi_1\xi_2 + c_n \xi_2^2 +d_{1n} \xi_1 + d_{2n} \xi_2 + e_n,
\end{equation}
with
$$
a_n = 2\pi^{-2}\int_0^1 \sin^2(\pi F_n(t)) \, dt,
$$
$$
b_n= \pi^{-2} \int_0^1 \sin (\pi F(t))\, \sin(2\pi F_n(t))\, dt,
$$
$$
c_n= \frac{1}{2} \pi^{-2} \int_0^1 \sin^2(2\pi F_n(t)) \, dt,
$$
$$
d_{1n} = \sqrt{2} \pi^{-1} \int_0^1 \sin(\pi F_n(t))\, J_n(t)\, dt,\quad d_{2n} = \frac{1}{ \sqrt{2}} \pi^{-1} \int_0^1 \sin(2\pi F(t_n))\, J_n(t)\, dt,
$$
$$
e_n = \int_0^1 J_n^2(t)\, dt,
$$
where
$$
J_n(t)
 = \int_0^1 \Bigl(\sqrt{2}\sum_{k=3}^\infty \xi_k \frac{\sin(\pi k F_n(t))}{k\pi} + n^{1/2}(F_n(t) - t)\Bigr)^2\, dt.
 $$
 We can write
\begin{equation}\label{pl3}
\mathbf{P}(T^2(\xi_1,\xi_2,J_n)  < c) = \int \chi_{\{T_n(x,y,\omega) < c\}}\,dG_n(x,y|\omega)\, d \mu_n(\omega),
\end{equation}
where $G_n(x,y|\omega)$ is conditional p.m. of $\xi_1$, $\xi_2$ given $J_n(t)$ and $\mu_n$  is p.m. of $J_n$.

Thus, for the proof of Lemma \ref{lom10} it suffices to show that distribution functions
\begin{equation}\label{pl4}
H_n(c|\omega) = \int \chi_{\{T_n(x,y,\omega) < c\}}\,dG_n(x,y|\omega)
\end{equation}
have uniformly bounded densities $ h_n(c|\omega)$ w.r.t. Lebesgue measure.

Define matrix $R_n =\{u_{ijn}\}_{i,j=1}^2$ with $u_{11n}= a_n$, $u_{22n}= c_n$ and $u_{12n} = u_{21n} = b_n$. Denote $I$ the unit matrix.

The distribution function $H_n(c|\omega)$ has characteristic function
\begin{equation}\label{pl5}
c \,(\mbox{\rm det}(I  -2itR))^{-1/2} \exp\{it q(a_n,b_n,c_n,d_{1n},d_{2n},e_n,\omega)\},
\end{equation}
where $q(a_n,b_n,c_n,d_{1n},d_{2n},e_n,\omega)$ is some function.

This is characteristic function of quadratic form of two Gaussian independent r.v.'s.
Therefore, if $\mbox{\rm det} (R_n)  > c$ then the densities $h_n$ are uniformly bounded.

We have
\begin{equation}\label{pl6}
\begin{split}&
\mbox{\rm det}(R_n) = \int^1_0\int^1_0\,(\sin^2(\pi F_n(x))\,\sin^2(2\pi F_n(y))\\& -
\sin(\pi F_n(x))\,\sin(\pi F_n(y))\,\sin(2\pi F_n(x))\,\sin(2\pi F_n(y)))\,dx\,dy\\&
= 4 \int^1_0\int^1_0\,(\sin^2(\pi F_n(x))\,\sin^2(\pi F_n(y))\, (cos^2(\pi F_n(y))\\& - cos(\pi F_n(x))\,\cos(\pi F_n(y))\,dx\,dy.
\end{split}
\end{equation}
Note that if we replace $\cos(\pi F_n(x))$ and $\cos(\pi F_n(y))$ with $|\cos(\pi F_n(x))|$ and $|\cos(\pi F_n(y))|$ respectively the right-hand side of (\ref{pl6}) remains nonnegative. Since $1 + f_n(x) \ge \delta$ and   $\cos(\pi F_n(x))$ and $\cos(\pi F_n(y))$ have both positive and negative values we get $\mbox{\rm det} (R_n) > c(\delta) > 0$. This completes the proof of Lemma \ref{lom10}.

Denote $F_{1n} = F_{n\gamma}$, $\gamma > 0$.

Since $T$ is a norm, by Hungary construction (see Th. 3, Ch. 12,  section 1, Schorack and Wellner \cite{wel}) and by Lemma \ref{lom10} the proof of (\ref{uuu}) and (\ref{uu1}) is reduced to the proof of two following inequalities.
\begin{equation}\label{pl7}
\begin{split}&
|\mathbf{P}( T^2(b(F_n(t)) +  \sqrt{n}(F_n(t) - F_0(t))) > x_\alpha)\\& - \mathbf{P}( T^2(b(F_{1n}(t)) +  \sqrt{n}(F_{1n}(t) - F_0(t))) > x_\alpha)| < \epsilon
\end{split}
\end{equation}
and
\begin{equation}\label{pl8}
\mathbf{P}(T^2(b(F_n(t)-F_{1n}(t) + F_0(t)) +  \sqrt{n}(F_n(t) - F_{1n}(t))) < x_\alpha) > 1 - \alpha - \epsilon.
\end{equation}
Since $T$ is a norm, the proof of (\ref{pl7}) and (\ref{pl8}) is reduced to the proof that, for any $\delta_1 > 0$, there hold
\begin{equation}\label{pl9}
\mathbf{P}(|T(b(F_n(t))) - T(b(F_{1n}(t)))| > \delta_1) = o(1),
\end{equation}
\begin{equation}\label{pl10}
\mathbf{P}(|T(b(F_0(t)+ F_n(t) - F_{1n}(t))) - T(b(F_{0}(t)))| > \delta_1) = o(1),
\end{equation}
and
\begin{equation}\label{pl11}
n^{1/2} |T(F_n(t)) - T(F_{1n}(t))|  <\delta_n(\gamma),
\end{equation}
\begin{equation}\label{pl12}
n^{1/2}|T(F_0(t)+ F_n(t) - F_{1n}(t)) - T(F_{0}(t))|  <\delta_n(\gamma),
\end{equation}
where $\delta_n(\gamma) \to 0$ as $\gamma \to \infty$ and $n \to \infty$.

Note that
\begin{equation}\label{pl13}
|T(b(F_n(t))) - T(b(F_{1n}(t)))| \le T(b(F_n(t)) - b(F_{1n}(t)))
\end{equation}
and
\begin{equation}\label{pl14}
 |T(F_n(t)) - T(F_{1n}(t))| \le T (F_n(t) - F_{1n}(t)).
\end{equation}
We have
\begin{equation}\label{pl15}
\begin{split}&
\mathbf{E} T^2(b(F_n) - b(F_{1n})) = \int_0^1 \mathbf{E}(b(F(t)) - b(F_{1n}(t)))^2 \, dt\\&
= \int_0^1 ((F_n(t) - \min(F_{n}(t), F_{1n}(t)) +
(F_{1n}(t) - \min(F_{n}(t), F_{1n}(t))\\& - (F_{n}(t)- F_{1n}(t))^2 \, dt \le C \max_{0<t<1}|F_n(t) - F_{1n}(t)|.
\end{split}
\end{equation}
In section 2 we point out Fourier coefficients $\eta_{nj}$ of functions $f_{n\gamma}$. In particular $\eta_{nj} = \theta_{nj}$ for $j < k_n$ with $k_n= [l_n]$ satisfying the equation
$$
l_n^{2s} \sum_{j=k_n}^\infty \theta_{nj}^2 = \gamma^2.
$$
Since
$$ cn^{-2r} \le \sum_{j=k_n}^\infty \theta_{nj}^2 \le C n^{-2r},$$ we get $\gamma^2 k_n^{-2s} < C n^{-2r}$. This implies
\begin{equation}\label{pl17}
k_n > (\gamma^2/C)^{1-2r} n^{\frac{1}{2}-r}.
\end{equation}
Therefore
\begin{equation}\label{pl18}
\begin{split}&
\max_{0<t<1}|F_n(t) - F_{1n}(t)| \le C \sum_{j=k_n}^\infty \frac{|\theta_{nj}|}{j}\\& \le C\left(\sum_{j=k_n}^\infty \theta_{nj}^2\right)^{1/2}  \left( \sum_{j=k_n}^\infty j^{-2}\right)^{1/2} \le Cn^{-r} k_n^{-1/2} \le C \gamma^{2r-1}.
\end{split}
\end{equation}
By (\ref{pl13}),(\ref{pl15}) and (\ref{pl18}), we get (\ref{pl9}).

We have
\begin{equation}\label{pl16}
\begin{split}&
 nT^2(F_n(t) - F_{1n}(t)) \le C \sum_{j=k_n}^\infty \frac{\theta_{nj}^2}{j^2} \\& \le C n k_n^{-2} \sum_{j=k_n}^\infty \theta_{nj}^2 \le  C n k_n^{-2} n^{-2r} \le C \gamma^{8r-4}.
 \end{split}
 \end{equation}
By (\ref{pl14}) and (\ref{pl16}),  we get (\ref{pl11}).

By Lemma \ref{lom10} and (\ref{pl9}), (\ref{pl11}), we  get (\ref{pl7}).

Proof of (\ref{pl10}) and (\ref{pl12}) is similar and is omitted.

\subsection{\bf Proof of Theorem \ref{t1}}
Fix $\delta, 0<\delta<1$. Denote $\kappa_j^2(\delta) = 0$ for $j > \delta^{-1}k_n$. Define $\kappa_j^2(\delta), 1 \le j < k_{n\delta}=\delta^{-1}k_n,$  the equations (\ref{i2}) and (\ref{i3}) with $P_0$ and $\rho_\epsilon$ replaced with $P_0(1-\delta)$ and $\rho_n(1 + \delta)$ respectively.
Similarly to \cite{er90}, we find Bayes test for a priori distribution $\theta_j = \eta_j=\eta_j(\delta), 1 \le j < \infty,$ with Gaussian independent random variables $\eta_j, E\eta_j =0, E\eta_j^2 = \kappa_j^2(\delta)$, and show that these tests are asymptotically minimax for some $\delta=\delta_n \to 0$ as $ n \to \infty$.
\begin{lemma}\label{l3} For any $\delta,  0 < \delta <1,$  there holds
\begin{equation}\label{i7}
\mathbf{P}(\eta(\delta) = \{\eta_j(\delta)\}_{j=1}^\infty  \in V_n) = 1 + o(1)
\end{equation}
as $n \to \infty$.
\end{lemma}
Denote
$$  A_{n,\delta} = \sigma^{-4}n^{2}\sum_{j=1}^\infty \kappa_j^4(\delta).
$$
By straightforward calculations, we get
\begin{equation}\label{i8}
\lim_{\delta \to 0}\lim_{n \to \infty} A_n A_n^{-1}(\delta) =1.
\end{equation}
Denote $\gamma_j^2(\delta) = \kappa_j^2(\delta)(n^{-1}\sigma^2 + \kappa_j^2(\delta))^{-1}$.

By Neymann-Pearson Lemma, Bayes critical region is defined the inequality
\begin{equation}\label{i9}
\begin{split}&
C_1 < \prod_{j=1}^{k_{n\delta}}(2\pi)^{-1/2}\kappa_j^{-1}(\delta) \int \exp\Bigl\{- \sum_{j=1}^{k_{n\delta}}(2\gamma_j^2(\delta))^{-1}(u_j- \gamma_j^2(\delta)y_j)^2\Bigr\}du \exp\{-T_{n\delta}(y)\}\\&
= C\exp\{-T_{n\delta}(y)\}(1+o(1))
\end{split}
\end{equation}
where
$$
T_{n\delta}(y) = n\sigma^{-2}\sum_{j=1}^\infty \gamma_j^2(\delta)y_j^2.
$$
Define critical region
$$
S_{n\delta}=\{y : \, R_{n\delta}(y) = (T_{n\delta}(y) -C_{n\delta})(2A_n(\delta))^{-1/2} > x_\alpha\}
$$
with
$$
C_{n\delta} = \mathbf{E}_0 T_{n\delta}(y) = \sigma^{-2}n\sum_{j=1}^\infty \gamma_j^2(\delta).
$$
Denote $L_{n\delta}$ the tests with critical regions $S_{n\delta}$.

Denote $\gamma_j^2 = \kappa_j^2(n^{-1} \sigma^2 + \kappa_j^2)^{-1}, 1 \le j < \infty$
 Define test statistics $T_{n}, R_n$,  critical regions $S_n$ and constants $C_n$  by the same way as test statistics $T_{n\delta}, R_{n\delta}$, critical regions $S_{n\delta}$ and constants $C_{n,\delta}$ respectively with $\gamma_j^2(\delta)$ replaced with $\gamma_j^2$ respectively.
Denote $L_n$ the test  having critical region $S_n$.
\begin{lemma}\label{l1} Let $H_0$ hold.  Then the distributions of tests statistics $R^a_n(y)$ and $R_n(y)$ converge to the standard normal distribution.

For any family $\theta_n =\{\theta_{jn}\} \in \Im_n$ there holds
\begin{equation}\label{i5a}
\mathbf{P}_{\theta_n}\Bigl(\Bigl(T_n^a(y) - C_n - \sigma^{-4}n^{2}\sum_{j=1}^\infty \kappa_j^2\theta_{jn}^2\Bigr)(2A_n)^{-1/2} < x_\alpha\Bigr) = \Phi(x_\alpha)(1+o(1))
\end{equation}
and
\begin{equation}\label{i6a}
\mathbf{P}_{\theta_n}\Bigl(\Bigl(T_n(y) - C_n - \sigma^{-4}n^2\sum_{j=1}^\infty \kappa_j^2\theta_{jn}^2\Bigr)(2A_n)^{-1/2} < x_\alpha\Bigr) = \Phi(x_\alpha)(1+o(1))
\end{equation}
as $n \to \infty$.
\end{lemma}
Hence we get the following Lemma.
\begin{lemma}\label{l2} There holds
\begin{equation}\label{i6}
\beta(L_n,V_n) = \beta(L^a_n,V_n )(1+ o(1))
\end{equation}

as $n \to \infty$.
\end{lemma}
\begin{lemma}\label{lx} Let $H_0$ hold.  Then the distributions of tests statistics $(T_{n\delta}(y) - C_{n\delta})(2A_n)^{-1/2}$ converge to the standard normal distribution.

There holds
\begin{equation}\label{i9a}
\mathbf{P}_{\eta(\delta)}((T_{n\delta}(y) - C_{n\delta} - A_{n\delta})(2A_{n\delta})^{-1/2} < x_\alpha) = \Phi(x_\alpha)(1+o(1))
\end{equation}
as $n \to \infty$.
\end{lemma}
\begin{lemma}\label{l4} There holds
\begin{equation}\label{i10}
\lim_{\delta\to 0}\lim_{n \to \infty} \mathbf{E}_{\eta(\delta)} \beta_{\eta(\delta)} (L_{n\delta})  =
\lim_{n \to \infty} \mathbf{E}_{\eta_0} \beta_{\eta_0} (L_{n})
\end{equation}
where $\eta_0 = \{\eta_{0j}\}_{j=1}^ \infty$ and $\eta_{0j}$ are i.i.d. Gaussian random variables, $\mathbf{E}[\eta_{0j}] =0$, $\mathbf{E} [\eta_{0j}^2] = \kappa_j^2, 1 \le j < \infty$.
\end{lemma}
Define Bayes a priori distribution  $\mathbf{P}_y$ as a conditional distribution of $\eta$ given $\eta \in V_n$. Denote $K_n = K_{n\delta}$ Bayes test with Bayes a priori distribution $P_y$. Denote $W_n$ critical region of $K_{n\delta}$.

For any sets $A$ and $B$ denote $A \triangle B = (A\setminus B)\cup(B\setminus A)$.
\begin{lemma}\label{l5} There holds
\begin{equation}\label{i11}
\lim_{\delta\to 0}\lim_{n \to \infty} \int_{V_n} \mathbf{P}_\theta(S_{n\delta}\triangle V_{n\delta}) d\,\mathbf{P}_y = 0
\end{equation}
and
\begin{equation}\label{i12}
\lim_{\delta\to 0}\lim_{n \to \infty} \mathbf{P}_0(S_{n\delta}\triangle V_{n\delta}) = 0.
\end{equation}
\end{lemma}
In the proof of Lemma \ref{l5} we show that the integrals in the right hand-side of (\ref{i9}) with integration domain $V_n$ converge to one in probability as $n \to \infty$. This statement is proved both for hypothesis and Bayes alternative (see \cite{er90}).

Lemmas \ref{l3}-\ref{l5} implies that, if $\alpha(K_n) = \alpha(L_n)$, then
\begin{equation}\label{i13}
\int_{V_n} \beta_\theta(K_n)\, d\mathbf{P}_y = \int_{V_n} \beta_\theta(L_n)\, d\mathbf{P}_y(1 + o(1)) = \int \beta_{\eta_0} (L_n)\, d\mathbf{P}_{\eta_0} (1+o(1)).
\end{equation}
\begin{lemma}\label{l6} There holds
\begin{equation}\label{i14}
\mathbf{E}_{\eta_0} \beta_{\eta_0}(L_n)  = \beta_n(L_n)(1+o(1)).
\end{equation}
\end{lemma}
Lemmas \ref{l1}, \ref{l4},  (\ref{i8}), (\ref{i13}) and Lemma \ref{l6}, imply Theorem \ref{t1}.
\subsection{Proof of Lemmas}
Proofs of Lemmas \ref{l1},  \ref{l2} and \ref{l4} are akin to the proofs of similar statements in \cite{er90} and are omitted.

\noindent{\sl Proof of Lemma \ref{l3}}. By straightforward calculations, we get
\begin{equation}\label{p1b}
\sum_{j=1}^\infty \mathbf{E}\eta_j^2(\delta) \ge \rho_\epsilon(1+\delta/2)
\end{equation}
and
\begin{equation}\label{p2}
\mathbf{Var}\Bigl(\sum_{j=1}^\infty \eta_j^2(\delta)\Bigr) < Cn^2A_n \asymp  \rho_n^2 k_n^{-1}.
\end{equation}
Hence, by Chebyshev inequality, we get
\begin{equation}\label{p3}
\mathbf{P}\Bigl(\sum_{j=1}^\infty \eta_j^2(\delta)> \rho_n\Bigr) = 1 +o(1)
\end{equation}
as $n \to \infty$.
It remains to estimate
\begin{equation}\label{e5}
\mathbf{P}_\mu(\eta \notin B^s_{2\infty}(P_0)) = \mathbf{P}\Bigl(\max_{l_1 \le i \le l_2} i^{2s} \sum_{j = i}^{l_2} \eta_j^2-P_0(1-\delta_1/2\Bigr) > P_0\delta_1/2) \le \sum_{i=l_1}^{l_2} J_i
\end{equation}
with
$$
J_i = \mathbf{P}\Bigl( i^{2s} \sum_{j = i}^{l_2} \eta_j^2-P_0(1-\delta_1/2)> P_0\delta_1/2\Bigr)
$$
To estimate $J_i$ we implement the following Proposition (see \cite{hsu}).
\begin{proposition}\label{p1} Let $\xi = \{\xi_i\}_{i=1}^l$ be Gaussian random vector with i.i.d.r.v.'s $\xi_i$, $\mathbf{E}[\xi_i] = 0$, $\mathbf{E}[\xi_i^2]=1$. Let $A\in R^l\times R^l$ and $\Sigma = A^T A$. Then
\begin{equation}\label{e6}
\mathbf{P}(||A\xi||^2 > \mbox{\rm tr}(\Sigma) + 2\sqrt{\mbox{\rm tr}(\Sigma^2)t} + 2 \|\Sigma\|t) \le \exp\{-t\}.
\end{equation}
\end{proposition}
We put $\Sigma_i= \{\sigma_{lj}\}_{l,j=i}^{k_{\epsilon\delta}}$ with $\sigma_{jj} = j^{-2s-1}i^{2s}\frac{P_0-\delta}{2s}$ and $\sigma_{lj} =0$ if $l\ne j$.

Let $i \le k_n$. Then
\begin{equation}\label{p4}
\mbox{\rm tr}(\Sigma_i^2) = i^{4s}\sum_{j=i}^\infty \kappa_j^4(\delta) < i^{4s}((k_n -i)\kappa^4(\delta) + k_n^{-4s-1}P_0) < Ck_n^{-1}.
\end{equation}
and
\begin{equation}\label{p5}
\|\Sigma_i\|  \le i^{2s}\kappa^2 < Ck_n^{-1}.
\end{equation}
Therefore
\begin{equation}\label{p6}
 2\sqrt{\mbox{\rm tr}(\Sigma_i^2)t} + 2 \|\Sigma_i\|t \le C(\sqrt{k_n^{-1}t} + k_n^{-1}t)
\end{equation}
Hence, putting $t =k_n^{1/2}$, by Proposition \ref{p1}, we get
\begin{equation}\label{p7}
\sum_{i=1}^{k_n} J_i \le Ck_n \exp\{-Ck_n^{1/2}\}.
\end{equation}
Let $i \ge k_n$.  Then
\begin{equation}\label{p8}
\mbox{tr}(\Sigma_i^2)  < Ci^{-1}, \quad\mbox{and} \quad ||\Sigma_i|| \le Ci^{-1}
\end{equation}
Hence, putting $t =i^{1/2}$, by Proposition \ref{p1}, we get
\begin{equation}\label{p9}
\sum_{i=k_n+1}^{k_{n\delta}} J_i \le \sum_{i=k_n+1}^{k_{n\delta}} \exp\{-Ci^{1/2}\} < \exp\{-C_1k_n^{1/2}\}.
\end{equation}
Now (\ref{e5}), (\ref{p7}), (\ref{p9}) together implies Lemma \ref{l3}.

\noindent{\sl Proof of Lemma \ref{l5}}. By reasoning of the proof of Lemma 4 in \cite{er90}, Lemma \ref{l5} will be proved, if we show, that
\begin{equation}\label{p10}
\mathbf{P}\Bigl(\sum_{j=1}^\infty (\eta_j(\delta) + y_j\gamma_j(\delta)\sigma^{-1}n^{1/2})^2 > \rho_n\Bigr) =1 +o(1)
\end{equation}
and
\begin{equation}\label{p11}
\mathbf{P}\Bigl(\sup_i i^{2s}\sum_{j=i}^\infty (\eta_j(\delta) + y_j\gamma_j(\delta)\sigma^{-1}n^{1/2})^2 > \rho_n\Bigr) =1 +o(1)
\end{equation}
where $y_j, 1 \le j < \infty$  are distributed by hypothesis  or Bayes alternative.

We prove only (\ref{p11}) in the case of Bayes alternative. In other cases the reasoning are similar.

We have
\begin{equation}\label{p12}
\begin{split}&
i^{2s}\sum_{j=i}^\infty (\eta_j(\delta) + y_j\gamma_j(\delta)\sigma^{-1}n^{1/2})^2=i^{2s}\sum_{j=i}^\infty \eta_j^2(\delta)\\& + i^{2s}\sum_{j=i}^\infty \eta_j(\delta)y_j\gamma_j(\delta)\sigma^{-1}n^{1/2} + i^{2s}\sum_{j=i}^\infty y_j^2\gamma_j^2(\delta)\sigma^{-2}n = J_{1i} + J_{2i}+ J_{3i}.
\end{split}
\end{equation}
The required probability for $J_{1n}$ is provided   Lemma \ref{l3}.

We have
\begin{equation}\label{p13}
J_{2i} \le J_{1i}^{1/2}J_{3i}^{1/2}.
\end{equation}
Thus it remains to show that, for any $C$,
\begin{equation}\label{p14}
\mathbf{P}_{\eta(\delta)}\Bigl(\sup_i i^{2s}\sum_{j=i}^\infty y_j^2\gamma_j^4(\delta)\sigma^{-2}n > C\delta\Bigr)   =   o(1)
\end{equation}
as $n \to \infty$.

Note that $y_j = \zeta_j + \sigma n^{-1/2}\xi_j$ where $\zeta_j, y_j, 1 \le j <\infty$ are i.i.d. Gaussian random variables, $\mathbf{E}\zeta_j =0, E\zeta_j^2 = \kappa_j^2(\delta), \mathbf{E} \xi_j =0, \mathbf{E} \xi_j^2 =1$.

Hence, we have
\begin{equation}\label{p15}
\begin{split}&
\sigma^{-2}n\sum_{j=i}^\infty y_j^2\gamma_j^4(\delta) =  \sigma^{-2}n\sum_{j=i}^\infty \gamma_j^4(\delta)\zeta_j^2  + \sigma^{-1}n^{1/2}\sum_{j=i}^\infty \gamma_j^4(\delta)\zeta_j\xi_j \\&+ \sum_{j=i}^\infty \gamma_j^4(\delta)\xi_j^2 = I_{1i}+ I_{2i}+ I_{3i}.
\end{split}
\end{equation}
Since $ n\gamma_j^2 = o(1)$, the estimates for probability of $i^{2s}I_{1i}$ are evident. It suffices to follow the estimates of (\ref{e5}). We have $I_{2i} \le I_{1i}^{1/2}I_{3i}^{1/2}$. Thus it remains to show that, for any $C$
\begin{equation}\label{p16}
\mathbf{P}_{\eta(\delta)}\Bigl( \sup_i i^{2s}\sum_{j=i}^\infty \gamma_j^4(\delta)\xi_j^2 > \delta/C\Bigr) =o(1)
\end{equation}
as $n \to \infty$.
Since $\gamma_j^2 = \kappa_j^2(1 + o(1)) = o(1)$, this estimate is also follows from estimates (\ref{e5}).

\noindent{\sl Proof of Lemma \ref{l6}}. By Lemmas \ref{l1}, \ref{l2} and \ref{l4}, it suffices to show that
\begin{equation}\label{p17}
\inf_{\theta \in V_n} \sum_{j=1}^\infty \kappa_j^2\theta_j^2 =  \sum_{j=1}^\infty \kappa_j^4.
\end{equation}
Denote $u_k = k^{2s}\sum_{j=k}^\infty \theta_j^2$. Note that $u_k \le P_0$.

Then $ \theta_j^2 = u_j j^{-2s} - u_{j+1}(j+1)^{-2s}$. Hence we have
\begin{equation}\label{p18}
\begin{split}&
A_n(\theta) = \sum_{j=1}^\infty \kappa_j^2\theta_j^2 = \kappa^2 \sum_{j=1}^{k_n} \theta_j^2 + \sum_{j=k_n}^\infty\kappa_j^2(u_j j^{-2s} - u_{j+1}(j+1)^{-2s})\\& =  \kappa^2 \sum_{j=1}^{k_n} \theta_j^2 + \kappa^2u_{k_n} k_n^{-2s} +  2s P_0 \sum_{j=k_n+1}^\infty u_j (j^{-4s-1} - (j-1)^{-2s-1}j^{-2s})\\& = \kappa^2 \rho_n + 2s P_0 \sum_{j=k_n+1}^\infty u_j (j^{-4s-1} - (j-1)^{-2s-1}j^{-2s}).
\end{split}
\end{equation}
Since $j^{-4s-1} - (j-1)^{-2s-1}j^{-2s}$ is negative, then $\inf A(\theta)$ is attained for $u_j = P_0$ and therefore $\theta_j^2 = \kappa_j^2$ for $j> k_\epsilon$.

Thus the problem is reduced  to the solution of the  following problem
\begin{equation}\label{p19}
\kappa^2\inf_{\theta_j}  \sum_{j=1}^{k_n} \theta_j^2 + \sum_{j=k_n+1}^\infty \kappa_j^4
\end{equation}
if
$$
\sum_{j=1}^{k_n} \theta_j^2 + \sum_{j=k_n+1}^\infty \kappa_j^2 =\rho_n
$$
and
$$
k_n^{2s}\sum_{j=k_n}^\infty \theta_j^2 < P_0, \quad 1 \le j < \infty,
$$
with $\theta_j^2 = \kappa_j^2$ for $j \ge k_n$.

It is easy to see that this infimum is attained if $\theta_j^2 =\kappa_j^2 = \kappa^2$ for $j \le k_n$.

\end{document}